# BETA-COALESCENTS AND CONTINUOUS STABLE RANDOM TREES

By Julien Berestycki, Nathanaël Berestycki and Jason Schweinsberg[1]

*Université de Provence, University of British Columbia and U.C. San Diego*

Coalescents with multiple collisions, also known as $\Lambda$-coalescents, were introduced by Pitman and Sagitov in 1999. These processes describe the evolution of particles that undergo stochastic coagulation in such a way that several blocks can merge at the same time to form a single block. In the case that the measure $\Lambda$ is the Beta$(2-\alpha, \alpha)$ distribution, they are also known to describe the genealogies of large populations where a single individual can produce a large number of offspring. Here, we use a recent result of Birkner et al. to prove that Beta-coalescents can be embedded in continuous stable random trees, about which much is known due to the recent progress of Duquesne and Le Gall. Our proof is based on a construction of the Donnelly–Kurtz lookdown process using continuous random trees, which is of independent interest. This produces a number of results concerning the small-time behavior of Beta-coalescents. Most notably, we recover an almost sure limit theorem of the present authors for the number of blocks at small times and give the multifractal spectrum corresponding to the emergence of blocks with atypical size. Also, we are able to find exact asymptotics for sampling formulae corresponding to the site frequency spectrum and the allele frequency spectrum associated with mutations in the context of population genetics.

**1. Introduction and preliminaries.** Consider the following simple population model. Assume that the size of the population stays constant, equal to a fixed integer $n \geq 1$, where individuals are numbered $1, \ldots, n$. In this

Received February 2006; revised July 2006.

[1]Supported in part by NSF Grant DMS-05-04882.

population, each individual reproduces at rate $(n-1)/2$. When individual $i$ reproduces, she gives birth to two children. One of them is again called individual $i$ and the other replaces individual $j$ for a randomly chosen label $j \neq i$ with $1 \leq j \leq n$. If $t > 0$ is a fixed time, we may define an *ancestral partition* $(\Pi_s^t, 0 \leq s \leq t)$ for this population model by saying that $i$ and $j$ are in the same block of $\Pi_s^t$ if and only if the corresponding individuals at time $t$ have the same ancestor at time $t-s$. It is elementary to check that the dynamics of the process $(\Pi_s^t, 0 \leq s \leq t)$ are governed by the rules of a process called *Kingman's coalescent*. This is a Markov process characterized by the fact that the only transitions are those where pairs of blocks merge and any given pair of blocks merges at rate 1 independently of everything else. In fact, even for more realistic population models, it is often the case that the genealogy of a small sample of a population may be effectively described by Kingman's coalescent; the introduction of this tool by Kingman [33, 34] was a major development in population genetics. One of the great advantages of this theory is that it is well adapted to the statistical analysis of molecular population samples since, for instance, in this framework, one can deal with a population sample rather than the population as a whole. Moreover, molecular and genetic data convey much information about ancestral relationships in a population sample. Much background material on the use of coalescent models in the field of population genetics can be found in the recent book [29] or in the review paper [25].

However, recent work (see, e.g., [21, 40, 49, 51]) has shown that Kingman's coalescent is not very well suited when we deal with populations where individuals may give birth to a large number of offspring or when we consider the genealogy of a population affected by repeated beneficial mutations [20]. In these cases, it is more appropriate to model the merging of ancestral lines by coalescent processes that allow multiple collisions, that is, several blocks may merge at once, although only one of those events may occur at a given time. These processes, called $\Lambda$-*coalescents*, have been introduced and studied by Pitman [45] and Sagitov [49]. As shown by Pitman [45], they are Markov processes in which any given number of blocks may merge at once and are characterized by a finite measure $\Lambda$ on $[0, 1]$. The $\Lambda$-coalescent has the property that whenever the process has $b$ blocks, any given $k$-tuple of blocks merges at a rate given by

$$\lambda_{b,k} = \int_0^1 x^{k-2}(1-x)^{b-k}\Lambda(dx);$$

see the next section for a more precise definition. For instance, Schweinsberg [51] showed that $\Lambda$-coalescents arise as the rescaled genealogies of some population models where individual offspring distributions have infinite variance. More precisely, let $1 < \alpha < 2$ and let $X$ be a random variable such that $P(X > k) \sim Ck^{-\alpha}$ for some $C > 0$. Consider the following population model.

As before, the size of the population is kept constant, equal to $n$. The model is formulated in discrete time. At each generation, each individual produces a random number of offspring (distributed like $X$) independently of other individuals and of the past. Then $n$ of them are randomly chosen to survive and the others are discarded. One of the main results of [51] is that the ancestral partitions, suitably rescaled, converge to the Beta$(2-\alpha, \alpha)$-coalescent, that is, a $\Lambda$-coalescent such that the measure $\Lambda$ is the Beta$(2-\alpha, \alpha)$ distribution.

This connection with population genetics has served both as a motivation for studying these processes and also as a source of inspiration for a rich theory that is only now starting to emerge, starting with the series of seminal papers by Bertoin and Le Gall [10, 11, 12, 13]. In these papers, $\Lambda$-coalescents are obtained as duals of measure-valued processes called *generalized Fleming–Viot processes*. In simple cases (viz., the cases of quadratic branching and stable branching mechanisms), these processes describe the composition of a population $(Z_t, t \geq 0)$ undergoing continuous branching (i.e., $Z$ is a continuous-state branching process, or CSBP for short; definitions will be given below). This stream of ideas has led Birkner et al. [14] to prove that one can obtain Beta-coalescents by suitably time-changing the ancestral partitions associated with the genealogy of $(Z_t, t \geq 0)$. In this continuous context, it is technically nontrivial to make rigorous sense of the notion of genealogy, but this is achieved through the use of a process called the (*modified*) *lookdown process* associated with $(Z_t, t \geq 0)$, a powerful tool introduced by Donnelly and Kurtz [16].

In parallel, it has been known for some time that CSBPs can be viewed as local time processes of a process $(H_t, t \leq T_r)$ called the *height process*, in a way that is analogous to the classical theorem of Ray and Knight for Brownian motion relating the Feller diffusion, the solution of

$$dZ_t = \sqrt{Z_t}\, dW_t,$$

where $(W_t)_{t \geq 0}$ is Brownian motion, to the local times of a reflecting Brownian motion. This connection has been formalized by Le Gall and Le Jan [36]. The height process itself encodes a continuous random tree, analogous to the Brownian tree of Aldous [1, 2], and can be viewed as the scaling limit of suitably normalized Galton–Watson trees. A careful exposition of this rich theory can be found in [17].

In this paper, we have two main goals. The first is to describe another way of thinking about the genealogy of a Beta-coalescent. This is achieved by embedding a Beta-coalescent into a continuous random tree with stable branching mechanism. To prove this result, we show that one can obtain the Donnelly–Kurtz lookdown process from a continuous random tree in a very simple fashion. This is valid for a general (sub)critical branching mechanism and is of independent interest. From this, and careful analysis, it follows that the coalescent tree can be thought of as what is perhaps the simplest genealogical model: a Galton–Watson tree with a continuous time parameter. Our second goal is to use this connection to discuss results about the small-time behavior of Beta-coalescents and related processes. This study was initiated in [8] without the help of continuous random trees. In particular, we apply these ideas to a problem of interest in population genetics.

*Organization of the paper.* After recalling the necessary definitions and results about coalescent processes, CSBPs and continuous random trees in Section 2, we state our results in Section 3. In Section 4, we explain our construction of the Donnelly–Kurtz lookdown process from a continuous random tree. In Section 5, we prove our results related to the small-time behavior of Beta-coalescents, giving asymptotics for the number of blocks and the multifractal spectrum. Finally, results concerning biological applications are proved in Section 6.

## 2. Preliminaries.

2.1. *The $\Lambda$-coalescent.* Let $\mathcal{P}_n$ denote the set of all partitions of the set $\{1,\ldots,n\}$ and $\mathcal{P}$ denote the set of all partitions of $\mathbb{N} = \{1, 2, \ldots\}$ (in this paper, it is always assumed that the set $\mathbb{N}$ does not contain 0). It turns out that the simplest way to define a coalescent process is by looking at a version of this process taking its values in the space $\mathcal{P}$. For all partitions $\pi \in \mathcal{P}$, let $R_n\pi$ be the restriction of $\pi$ to $\{1,\ldots,n\}$, meaning that $R_n\pi \in \mathcal{P}_n$ and that two integers $i$ and $j$ are in the same block of $R_n\pi$ if and only if they are in the same block of $\pi$. A $\Lambda$-coalescent (or a coalescent with multiple collisions) is a $\mathcal{P}$-valued Markov process $(\Pi(t), t \geq 0)$ such that, for all $n \in \mathbb{N}$, the process $(R_n\Pi(t), t \geq 0)$ is a $\mathcal{P}_n$-valued Markov chain with the property that whenever $R_n\Pi(t)$ has $b$ blocks, any particular $k$-tuple of blocks of this partition merges at a rate equal to $\lambda_{b,k}$, these being the only possible transitions. The rates $\lambda_{b,k}$ depend neither on $n$ nor on the numbers of integers in the $b$ blocks. Pitman [45] showed that the transition rates must satisfy

$$\lambda_{b,k} = \int_0^1 x^{k-2}(1-x)^{b-k}\Lambda(dx) \tag{1}$$

for some finite measure $\Lambda$ on $[0,1]$. The laws of the processes $R_n\Pi$ are consistent and this allows one to consider a process $\Pi$ such that the restriction $R_n\Pi$ has the above description. A coalescent process such that (1) holds for a particular measure $\Lambda$ is called the $\Lambda$-*coalescent*.

To better understand the role of the measure $\Lambda$, it is useful to have in mind the following Poissonian construction of a $\Lambda$-coalescent, also due to Pitman [45]. Suppose $\Lambda$ does not put any mass on $\{0\}$. Let $(t_i, x_i)_{i \in I}$ be the atoms of a Poisson point process on $\mathbb{R}^+ \times [0,1]$ with intensity measure $dt \otimes x^{-2}\Lambda(dx)$. Observe that although $\Lambda$ is a finite measure, $x^{-2}\Lambda(dx)$ is not finite in general, but only sigma-finite. Hence, $(t_i, x_i)_{i \in I}$ may have countably many atoms on any time interval $[t_1, t_2]$, so in order to make rigorous sense of the following description, one should again work with restrictions to $\{1, \ldots, n\}$. The coalescent only evolves at times $t$ such that $t = t_i$ for some $i \in I$. For each cluster present at time $t_i^-$, we flip an independent coin with probability of heads $x_i$, where $(t_i, x_i)$ is the corresponding atom of the point process. We merge all the clusters for which the coin came up heads and do nothing with the other clusters. Hence, we see that in a $\Lambda$-coalescent where $\Lambda$ has no mass at 0, $x^{-2}\Lambda(dx)$ is the rate at which a proportion $x$ of the blocks merge (such an event is generally called an $x$-*merger*). On the other hand, when $\Lambda$ is a unit mass at zero, each transition involves the merger of exactly two blocks and each such transition occurs at rate 1, so this is just Kingman's coalescent.

Kingman's theory of exchangeable partitions provides us with a way of looking at this process as taking its values in the space $\mathcal{S} = \{x_1 \geq x_2 \geq \cdots \geq 0, \sum_{i=1}^\infty x_i \leq 1\}$, which is perhaps a bit more intuitive since the notion of mass is apparent in this context. The resulting process is called the *ranked $\Lambda$-coalescent*. Briefly, partitions of $\mathbb{N}$ defined by the above procedure are exchangeable, so this implies that for each block of the partition, there exists a well-defined number called the *frequency* or *mass* of the block, which is the almost sure limiting proportion of integers in this block. Therefore, given a measure $\Lambda$ and a $\Lambda$-coalescent $\Pi = (\Pi_t, t \geq 0)$, one can define a process $X = (X(t), t > 0)$ with values in the space $\mathcal{S}$ by taking for each $t > 0$ the frequencies of $\Pi(t)$ ranked in decreasing order. When $\mathcal{S}$ is endowed with the topology that it inherits from $\ell^1$, the law at time $t$ of this process $Q_t$ defines a Markov semigroup with an *entrance law*: the process enters at time $0^+$ from a state called *dust*, that is, the largest frequency vanishes as $t \to 0^+$. These technical points are carefully explained in the original paper of Pitman [45], Theorem 8. The process $X$ is said to have *proper frequencies* if $\sum_{i=1}^\infty X_i(t) = 1$ for all $t > 0$. Pitman has shown that this is equivalent to $\int_0^1 x^{-1}\Lambda(dx) < \infty$. This is also equivalent to the fact that almost surely $\Pi(t)$ does not contain any singleton, or that all blocks are infinite. Another notion which plays an important role in this theory is that of *coming down from infinity*. Pitman [45] has shown that only two situations occur, depending on the measure $\Lambda$. Let $E$ be the event that for all $t > 0$, there are infinitely many blocks and let $F$ be the event that for all $t > 0$, there are only finitely many blocks. Then, if $\Lambda(\{1\}) = 0$, either $P(E) = 1$ or $P(F) = 1$. When $P(F) = 1$, the process $X$ or $\Pi$ is said to come down from infinity. For instance, Kingman's coalescent comes down from infinity, while if $\Lambda(dx) = dx$ is the uniform measure on $(0,1)$, the $\Lambda$ coalescent does not come down from infinity. This particular choice of $\Lambda$ corresponds to the so-called Bolthausen–Sznitman coalescent which first arose in connection with spin glasses [15]. For a necessary and sufficient condition on $\Lambda$ for coming down from infinity, see [13, 50] and the forthcoming [7]. Note also that a coalescent that comes down from infinity must have proper frequencies.

In this paper, we will be concerned with the one-parameter family of coalescent processes called Beta-coalescents. These are the $\Lambda$-coalescent process obtained when the measure $\Lambda$ is the Beta$(2-\alpha,\alpha)$ distribution with $1 < \alpha < 2$,

$$\Lambda(dx) = \frac{1}{\Gamma(2-\alpha)\Gamma(\alpha)} x^{1-\alpha}(1-x)^{\alpha-1}\,dx.$$

The reason we restrict our attention to $1 < \alpha < 2$ is that this corresponds to the case where the coalescent process comes down from infinity (a consequence of Schweinsberg's [50] criterion). When $\alpha = 1$, the Beta$(1, 1)$ distribution is simply the uniform distribution on $(0, 1)$, so this the Boltahusen–Sznitman coalescent, which stays infinite. When $\alpha \to 2$, it can be checked that the Beta$(2-\alpha,\alpha)$ distribution converges weakly to the unit mass at zero, so, formally, the case $\alpha = 2$ corresponds to Kingman's coalescent. This family of processes enjoys some remarkable properties, as can be seen from [14, 51] and results in the present work. This partly reflects the fact that the continuous-state branching processes with stable branching mechanism, with which they are associated (see below), enjoy some strong scale-invariance properties, just like Brownian motion.

2.2. *Continuous-state branching processes.* Continuous-state branching processes have been introduced and studied by, among others, Lamperti [35] and Grey [27]. They are Markov processes $(Z_t, t \geq 0)$ taking their values in $[0, \infty]$ and we think of $Z_t \geq 0$ as the size of a continuous population at time $t$. Continuous-state branching processes are the continuous analogues of Galton–Watson processes as well as their scaling limits. They are characterized by the following branching property: if $p_t(x,\cdot)$ denotes the transition probabilities of $Z$ started with $Z_0 = x$, then for all $x, y \in \mathbf{R}_+$,

$$(2) \qquad p_t(x+y,\cdot) = p_t(x,\cdot) * p_t(y,\cdot),$$

which means that the process started from $x+y$ individuals has the same law as the sum of a process started from $x$ and one started from $y$ independently. The interpretation of (2) is that if individuals live and reproduce independently, then a population started from $x+y$ individuals should evolve as the sum of two independent populations, one started with $x$ individuals and one with $y$ individuals. Lamperti [35] has shown that a continuous-state branching process is characterized by a function $\psi\colon[0,\infty)\to\mathbb{R}$ called the *branching mechanism*, such that for all $t\geq 0$, the Laplace transform of $Z_t$ satisfies

$$E[e^{-\lambda Z_t}|Z_0=a]=e^{-au_t(\lambda)}, \tag{3}$$

where the function $u_t(\lambda)$ solves the differential equation

$$\frac{\partial u_t(\lambda)}{\partial t}=-\psi(u_t(\lambda)), \qquad u_0(\lambda)=\lambda. \tag{4}$$

Moreover, the branching mechanism $\psi$ is the Laplace exponent of some spectrally positive Lévy process (i.e., Lévy process with no negative jumps). That is, there exists a measure $\nu$ on $(0,\infty)$ and some numbers $a\in\mathbb{R}$ and $b\geq 0$ such that for all $q\geq 0$,

$$\psi(q)=aq+bq^2+\int_0^\infty(e^{-qx}-1+qx\mathbf{1}_{\{x<1\}})\nu(dx) \tag{5}$$

and $\int_0^\infty(1\wedge x^2)\nu(dx)<\infty$. Furthermore, if $(Y_t,t\geq 0)$ is the Lévy process with Laplace exponent $\psi$, that is,

$$E[e^{-\lambda(Y_t-Y_0)}]=e^{t\psi(\lambda)},$$

then the distributions of $(Z_t,t\geq 0)$ and $(Y_t,t\geq 0)$ are related by a simple time-change (sometimes called the *Lamperti transform*). Let

$$U_t=\int_0^t \tilde{Y}_s^{-1}\,ds,$$

where $(\tilde{Y}_t,t\geq 0)$ is the process $(Y_t,t\geq 0)$ stopped when it first hits zero, and call $U_t^{-1}$ the inverse càdlàg of $U_t$. $(Y_{U_t^{-1}},t\geq 0)$ then has the same law as $Z$. We refer the reader to, for instance, [9] for more information about this.

When $\psi(q)=q^\alpha$ for some $\alpha\in(0,2]$, we say that the CSBP has a *stable* branching mechanism. When $\alpha=2$, the process $Z$ is Feller's diffusion and the Lévy process in Lamperti's transformation is standard Brownian motion. When $1<\alpha<2$, this branching mechanism arises by taking $a=b=0$ and

$$\nu(dx)=\frac{\alpha(\alpha-1)}{\Gamma(2-\alpha)}x^{-1-\alpha}\,dx$$

in (5). The Lévy process in Lamperti's transformation is an $\alpha$-stable Lévy process having the scaling property

$$(Y_{\lambda t},t\geq 0)=_d(\lambda^{1/\alpha}Y_t,t\geq 0).$$

2.3. *The height process and continuous random trees.* Le Gall and Le Jan [36] have introduced a new way of thinking about CSBPs, which was further carefully explored by Duquesne and Le Gall in [17]. It is inspired by the well-known result of Ray and Knight on the local times of Brownian motion and is related to the construction of the Brownian continuum random tree of Aldous [1, 2]. Recall that if $B$ is a reflecting Brownian motion, $(L_t^x, t \geq 0, x \geq 0)$ is a jointly continuous version of its local times and $(T_r, r \geq 0)$ is the càdlàg inverse of $L_t^0$, then for fixed $r > 0$, the process $(L_{T_r}^x, x \geq 0)$ is a Feller diffusion started with initial population $r$. Le Gall and Le Jan have introduced a process $(H_t, t \geq 0)$ which generalizes the Ray–Knight theorem to continuous branching process with (sub)critical branching mechanism.

More precisely, consider a Laplace exponent $\psi(q)$ and a $\psi$-CSBP $(Z_t, t \geq 0)$. We will assume that $\psi$ is subcritical, that is, a.s. there exists some time $0 < \tau < \infty$ such that $Z_\tau = 0$. Grey has shown that this is equivalent to the condition that the branching mechanism $\psi$ satisfies

$$\int_1^\infty \frac{dq}{\psi(q)} < \infty.$$

In particular, this is the case when $\psi(q) = q^2/2$ or when $\psi(q) = q^\alpha$ for $1 < \alpha < 2$. Lamperti [35] has shown that there exists a sequence of offspring distributions $\mu_n$ such that if we consider $(Z_k^n, k = 1, 2, \ldots)$, a discrete Galton–Watson process with offspring distribution $\mu_n$ and started with $n$ individuals, then $(n^{-1} Z_{\gamma_n t}^n, t \geq 0)$ converges in the sense of finite-dimensional distributions to $(Z_t, t \geq 0)$, where the $\gamma_n$ are suitable time-scaling constants. If we ask for finer limit theorems about the genealogy of $(Z_t, t \geq 0)$, then Duquesne and Le Gall have shown that the discrete height process $(H_k^n, k = 0, 1, \ldots)$, where $H_k^n$ is the generation of the $k$th individual, converges when suitably normalized to a process $(H_t, t \geq 0)$ called the *height process*. One may directly construct this process $(H_t, t \geq 0)$ from a Lévy process with Laplace exponent $\psi$. Thus, informally, the height process plays the same role as the depth-first search process on a discrete tree, but in a continuous setting. An important result of Duquesne and Le Gall [17] is that, even though $H$ is, in general, neither a semi-martingale nor a Markov process, that is, $H$ admits a local time process, that is, almost surely, there exists a jointly continuous process $(L_s^a, s \geq 0, a \geq 0)$ such that for all $t \geq 0$,

$$\lim_{\varepsilon \to 0} E\left[\sup_{0 \leq s \leq t} \left| \frac{1}{\varepsilon} \int_0^s \mathbf{1}_{\{a < H_r \leq a+\varepsilon\}} \, dr - L_s^a \right|\right] = 0.$$

They were also able to prove that the process $H$ has a continuous modification provided the branching mechanism is subcritical.

The importance of the process $H$ stems primarily from the generalized Ray–Knight theorem, which we now state (see [17, 36]). Let $T_r = \inf\{t >$

$0, L_t^0 > r\}$ be the inverse local time at 0. For all $t \geq 0$, define

$$Z_t = L_{T_r}^t. \tag{6}$$

Then $(Z_t, t \geq 0)$ is a $\psi$-CSBP started at $Z_0 = r$. If $\psi(q) = q^2/2$, then $(H_t, t \geq 0)$ has the law of a reflecting Brownian motion and $(Z_t, t \geq 0)$ is the Feller diffusion, as the classical Ray–Knight theorem states.

## 3. Main results.

3.1. *The Beta-coalescent in the continuous stable random tree.* Our first result is the embedding of a Beta$(2-\alpha, \alpha)$-coalescent for $1 < \alpha < 2$ in the tree coded by the $\alpha$-stable height process. Let $Z$ be an $\alpha$-stable CSBP obtained in the fashion of Duquesne and Le Gall from the height process $(H_t, 0 \leq t \leq T_r)$ associated with $\psi(q) = q^\alpha$ for a given $1 < \alpha < 2$, that is, $Z_t = L_{T_r}^t$. Consider, for all $t$, the random level

$$R_t = \alpha(\alpha - 1)\Gamma(\alpha) \int_0^t Z_s^{1-\alpha}\, ds \tag{7}$$

and let $R^{-1}(t) = \inf\{s : R_s > t\}$. It follows from [14] that $R^{-1}(t) < \infty$ a.s. for all $t$ and that $\lim_{t \to \infty} R^{-1}(t) = \zeta$, where $\zeta$ is the lifetime of the CSBP.

Let $(V_i, i = 1, 2, \ldots)$ be a sequence of variables in $(0, T_r)$ defined such that for all $i \in \mathbb{N}$, $V_i$ is the left endpoint of the $i$th highest excursion of the height process $H$ above the level $R^{-1}(t)$. Next, we define a process $(\Pi_s, 0 \leq s \leq t)$ which takes its values in the space $\mathcal{P}$ of partitions of $\mathbb{N}$ as follows:

$$i \overset{\Pi_s}{\sim} j \iff \left(\inf_{r \in [V_i, V_j]} H_r\right) > R^{-1}(t-s).$$

That is, $i$ and $j$ are in the same block of $\Pi_s$ if and only if $V_i$ and $V_j$ are in the same excursion of $H$ above level $R^{-1}(t-s)$.

THEOREM 1. *The process $(\Pi_s, 0 \leq s \leq t)$ is a Beta$(2-\alpha, \alpha)$-coalescent run for time $t$.*

Another way of looking at this result is to consider the ranked coalescent. Let $(X(s), 0 \leq s \leq t)$ be the process with values in $\mathcal{S}$ defined by the following procedure. For each $s \leq t$, $X(s)$ has as many nonzero coordinates as there are excursions of the height process above $R^{-1}(t-s)$ that reach the level $R^{-1}(t)$. To each such excursion we associate a mass given by the local time of that excursion at level $R^{-1}(t)$, normalized by $Z_{R^{-1}(t)}$ so that the sum is equal to 1. Then $X(s)$ is defined as the nonincreasing rearrangement of these masses.

COROLLARY 2. *$(X(s), 0 \leq s \leq t)$ has the same distribution as the ranked Beta$(2-\alpha, \alpha)$-coalescent run for time $t$.*

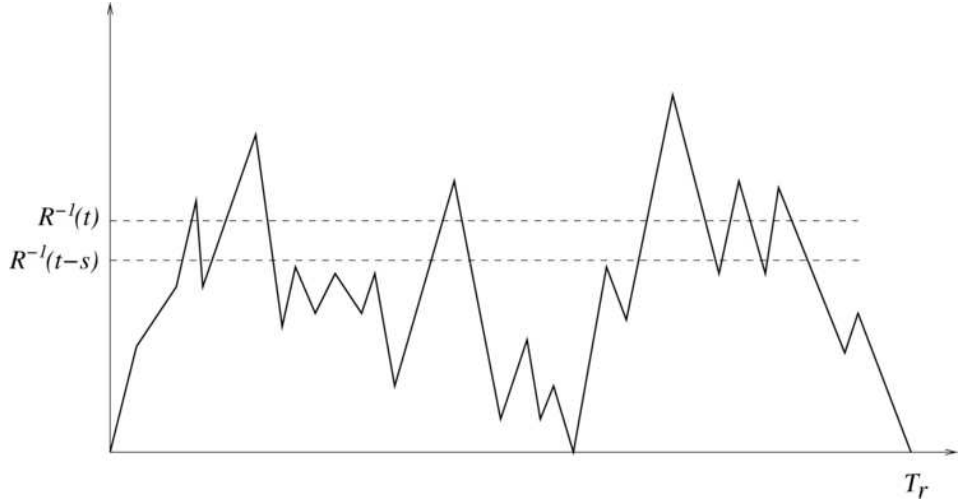

FIG. 1. *A Beta-coalescent is obtained by coalescing excursions of $(H_t, t \leq T_r)$ above $R^{-1}(t-s)$ that reach $R^{-1}(t)$. Thus each excursion corresponds to a block of the coalescent and its mass is given by its local time at level $R^{-1}(t)$.*

We picture the coalescent as the following process. As $s$ increases from 0 to $t$, the level $R^{-1}(t-s)$ decreases from $R^{-1}(t)$ to 0. The excursions of $H$ above level $R^{-1}(t-s)$ coalesce because if $s_1 < s_2$, then several excursions of $H$ above the level $R^{-1}(t-s_1)$ could be part of the same excursion of $H$ above the level $R^{-1}(t-s_2)$. This will happen, for example, if the excursion of $H$ above the level $R^{-1}(t-s_1)$ has a local minimum at the level $R^{-1}(t-s_2)$. Then, in the corresponding coalescent process, we observe a merging of masses at time $s_2$ corresponding to the fraction of local time at $R^{-1}(t)$ contained by each of those excursions.

REMARK 3. Recall the definition of an $\mathbb{R}$-tree associated with a nonnegative function $H$ defined on an interval $[0, T_r]$. If $d_H(u,v) = H(u) + H(v) - 2\inf_{u \leq t \leq v} H(t)$, then $d_H$ is a pseudodistance on $[0, T_r]$. Equipped with $d_H$, the quotient of $[0, T_r]$ by the relation $d_H(u,v) = 0$ is an $\mathbb{R}$-tree. For the function $(H_s, s \leq T_r)$, this gives a Poissonian collection of scaled stable trees joined at the root. In this context, the $V_i$ are certain vertices at distance $R^{-1}(t)$ from the root and the state of the coalescent at time $s$ can be described as the partition obtained by declaring $i \sim j$ if and only if their most recent common ancestor is at distance greater than $R^{-1}(t-s)$ from the root, that is, if $d_H(V_i, V_j) < 2(R^{-1}(t) - R^{-1}(t-s))$. In other words, if we define a new distance $d_\Pi^{(t)}$ on $\mathbb{N}$ by

$$d_\Pi^{(t)}(i,j) = \inf\{s : R^{-1}(t) - R^{-1}(t-s) = d_H(V_i, V_j)/2\},$$

then the classes of $\Pi_s$ are the balls of radius $s$ for the metric $d_\Pi^{(t)}$.

3.2. *Small-time behavior and multifractal spectrum.* We now use Theorem 1 to obtain several results about the small-time behavior of the Beta coalescents.

Let $N(t)$ be the number of blocks at time $t$ of a Beta-coalescent $\Pi(t)$. Our first application gives the almost sure limit behavior of $N(t)$ and has already been shown in [8] using methods based on the analysis of CSBP with stable branching mechanisms.

THEOREM 4.
$$\lim_{t \to 0} t^{1/(\alpha-1)} N(t) = (\alpha \Gamma(\alpha))^{1/(\alpha-1)} \qquad a.s.$$

For an exchangeable random partition, the number of blocks is related to the typical block size. For instance, suppose $\Pi$ is an exchangeable random partition and that $|\Pi|$ denotes the number of blocks of $\Pi$. Using equation (2.27) in [46], we see that if $X_1$ is the asymptotic frequency of the block of $\Pi$ containing 1, then $E(|\Pi|) = E(X_1^{-1})$. Hence, here, at least informally, we see that the frequency of the block which contains 1 at time $t$ must be of the order of $1/N(t) \propto t^{1/(\alpha-1)}$ (this result was proved rigorously in [8]). Put another way, this says that almost all of the fragments emerge from the original dust by growing like $t^{1/(\alpha-1)}$. We say that $1/(\alpha-1)$ is the *typical speed of emergence.*

However, some blocks clearly have a different behavior. Consider, for instance, the largest block and denote by $W(t)$ its frequency at time $t$. It was shown in [8], Proposition 1.6, that
$$(\alpha \Gamma(\alpha) \Gamma(2-\alpha))^{1/\alpha} t^{-1/\alpha} W(t) \to_d X \qquad \text{as } t \downarrow 0,$$
where $X$ has the Fréchet distribution of index $\alpha$. Hence, the size of the largest fragment is of the order of $t^{1/\alpha}$.

This suggests studying the existence of fragments that emerge with an atypical rate $\gamma \neq 1/(\alpha-1)$. To do so, it is convenient to consider a random metric space $(S, d)$ which completely encodes the coalescent $\Pi$ (this space was introduced by Evans [23] in the case of Kingman's coalescent). The space $(S, d)$ is the completion of the space $(\mathbb{N}, d)$, where $d(i, j)$ is the time at which the integers $i$ and $j$ coalesce. In particular, completing the space $\{1, 2, \ldots\}$ with respect to this distance in particular adds points that belong to blocks behaving atypically. In this framework, we are able to associate with each point $x \in S$ and each $t > 0$ a positive number $\eta(x, t)$ which is equal to the frequency of the block at time $t$ corresponding to $x$. (This is formally achieved by endowing $S$ with a mass measure $\eta$.) In this setting, we can reformulate the problem as follows: are there points $x \in S$ such that the block $B_x(t)$ that contains $x$ at time $t$ behaves as $t^\gamma$ when $t \to 0$ or, more formally, such that $\eta(x,t) \asymp t^\gamma$? [Here, $f(t) \asymp g(t)$ means that $\log f(t)/\log g(t) \to 1$.] Also, how many such points typically exist?

We define, for $\gamma \leq 1/(\alpha - 1)$,
$$S(\gamma) = \left\{ x \in S : \liminf_{t \to 0} \frac{\log(\eta(x,t))}{\log t} \leq \gamma \right\}$$
and, similarly, when $\gamma > 1/(\alpha - 1)$,
$$S(\gamma) = \left\{ x \in S : \limsup_{t \to 0} \frac{\log(\eta(x,t))}{\log t} \geq \gamma \right\}.$$

When $\gamma \leq 1/(\alpha - 1)$, $S(\gamma)$ is the set of points which correspond to large fragments. On the other hand, when $\gamma \geq 1/(\alpha - 1)$, $S(\gamma)$ is the set of points which correspond to small fragments. In the next result, we answer the question raised above by computing the Hausdorff dimension (with respect to the metric of $S$) of the set $S(\gamma)$.

THEOREM 5.

1. If $\frac{1}{\alpha} \leq \gamma < \frac{1}{\alpha-1}$, then
$$\dim_\mathcal{H} S(\gamma) = \gamma\alpha - 1.$$

    If $\gamma < 1/\alpha$, then $S(\gamma) = \varnothing$ a.s. but $S(1/\alpha) \neq \varnothing$ almost surely.

2. If $\frac{1}{\alpha-1} < \gamma \leq \frac{\alpha}{(\alpha-1)^2}$, then
$$\dim_\mathcal{H} S(\gamma) = \frac{\alpha}{\gamma(\alpha-1)^2} - 1.$$

    If $\gamma > \frac{\alpha}{(\alpha-1)^2}$, then $S(\gamma) = \varnothing$ a.s. but $S(\frac{\alpha}{(\alpha-1)^2}) \neq \varnothing$ almost surely.

REMARK 6. The maximal value of $\dim_\mathcal{H} S(\gamma)$ is obtained when $\gamma = 1/(\alpha - 1)$, in which case the dimension of $S(\gamma)$ is also equal to $1/(\alpha - 1)$. This was to be expected since this is the typical exponent for the size of a block. The value of the dimension then corresponds to the full dimension of the space $S$, as was proved in [8], Theorem 1.7.

REMARK 7. We recover part of Proposition 1.6 in [8] that the largest block has size of order $t^{1/\alpha}$ since this is the smallest $\gamma$ for which $S(\gamma) \neq \varnothing$. It may be a bit more surprising that there is such a thing as a notion of *smallest* block, whose size is of order $t^\gamma$, where $\gamma = \alpha/(\alpha - 1)^2$.

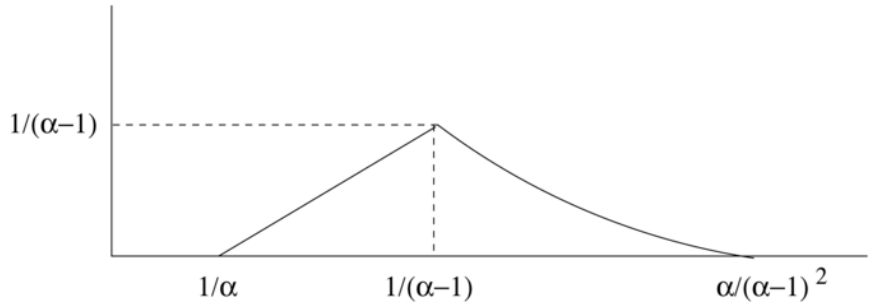

FIG. 2. *Multifractal spectrum map $\gamma \mapsto \dim_{\mathcal{H}} S(\gamma)$. The left-derivative at the critical point is $\alpha$ while the right-derivative is $-\alpha$.*

REMARK 8. This is reminiscent of the problem considered in [6], in which the long-time asymptotic behavior of homogeneous fragmentations was studied. More precisely, it was shown there that if $F(t)$ is a homogeneous fragmentation of the interval $(0,1)$ and $I_x(t)$ denotes the fragment that contains $x$ at time $t$, then there is a typical speed of fragmentation $v_0$, in the sense that if $U$ is uniform on $(0,1)$, then almost surely $|I_U(t)| \sim e^{-v_0 t}$. However, for $v \neq v_0$ in some range, the random set of exceptional points $S(v) := \{x \in (0,1) : |I_x(t)| \sim e^{-vt}\}$ is nonempty and has zero Lebesgue measure. The main result in [6] gives an explicit expression of the multifractal spectrum map $v \mapsto \dim_{\mathcal{H}}(S(v))$ where $\dim_{\mathcal{H}}(S)$ denotes the Hausdorff dimension of $S$. However, we emphasize that in [6], this Hausdorff dimension is computed with respect to the metric $\delta$ induced by the Lebesgue measure on $(0,1)$. In that case, the fact that the diameter of a block is equal to its mass plays a significant role. By contrast, here, we compute dimensions with respect to the metric $d$, which should rather be understood as a genealogical distance.

3.3. *Frequency spectra for mutation models.* We now describe a result concerning Beta-coalescents which has some applications to a question arising in population genetics. The question is concerned with the quantification of polymorphism in a sample of given size taken from a population. Suppose we sample $n$ individuals from a population at a certain time. Due to mutations, at a given locus, not all individuals in this sample will have the same allele. Moreover, mutations also affect different sites. We may ask several questions. In the sample of size $n$, how many different alleles should we observe at a given locus (site)? On how many sites should we expect to see different alleles? With which frequency should each of the different alleles be represented? As we will see, the answers to these questions depend heavily on the nature of the population, particularly on its reproduction mechanism, in addition to the mutation rate.

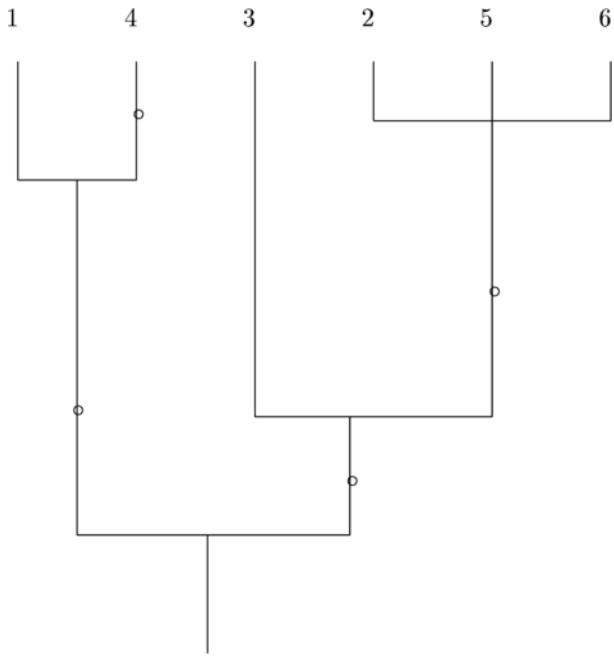

Fig. 3.   *In the infinite sites model, each mark stands for a mutation that affects a different locus. In this example, there are four families:* $\{1,4\}$, $\{4\}$, $\{2,3,5,6\}$ *and* $\{2,5,6\}$. *On the other hand, in the infinite alleles model, the allelic partition* $\Pi_\theta$ *also has four blocks:* $\{1\}$, $\{2,5,6\}$, $\{3\}$ *and* $\{4\}$.

To make the problem mathematically tractable, we will consider two simplified models. The rate at which mutations occur will always be assumed to be a positive number $\theta$, constant with time. In the first model, called the *infinite alleles model*, introduced by Kimura and Crow [32] in 1964, we study a given locus in the sample and assume that each mutation has resulted in a new allele. This means that the descendants of an individual affected by a mutation all carry the same allele except those later affected by another mutation. In the second model, called the *infinite sites model*, introduced by Kimura [31] in 1969, we look at the number of sites where we expect individuals to show polymorphism. In this model, we assume that each mutation occurs at a new site. In particular, if an individual is affected by a mutation, *all* the descendants of this individual carry this mutation. See Figure 3 for an illustration of these two models.

In the infinite alleles model, one can define the so-called *allelic partition*. That is, one may divide the sample into groups of individuals having the same allele at the observed locus. For a sample of size $n$, quantities of interest include the number of different groups (which we will also refer to as *clusters* or also sometimes *blocks*), $N(n)$, as well as typical sizes of groups: we denote by $N_k(n)$ the number of blocks in the allelic partition of size $k$. In the infinite sites model, one cannot define a partition of the sample because a given individual in the sample may have been affected by several mutations. However, we can still define $M(n)$ to be the total number of mutations and $M_k(n)$ to be the number of these mutations affecting exactly $k$ individuals in the sample. For example, in Figure 3, $N(n) = 4$, $N_1(n) = 3, N_2(n) = 0, N_3(n) = 1$, while $M(n) = 4$ and $M_1(n) = M_2(n) = M_3(n) = M_4(n) = 1$. The whole sequence $(M_1(n), \ldots, M_n(n))$ is called the *site frequency spectrum* and the sequence $(N_1(n), \ldots, N_n(n))$ is called the *allele frequency spectrum*.

A fundamental result in this domain is the celebrated Ewens sampling formula [24]. This result gives an explicit formula for the distribution of the allelic partition, under some standard assumptions on the reproduction mechanism of the population. The result is perhaps best explained through the theory of Kingman's coalescent. Based on this process, Kingman [33] was able to find a simpler proof of Ewens' sampling formula. Assume that the genealogy of the population may be described by the dynamics of Kingman's coalescent, that is, each pair of lineages coalesces at rate 1. Assuming the rate of mutations is $\theta/2$ along every lineage, the Ewens sampling formula states that the probability that the allelic partition has $a_i$ blocks of size $i$ for $i = 1, \ldots, n$ is

$$(8) \quad P(N_1(n) = a_1, \ldots, N_n(n) = a_n) = p(a_1, \ldots, a_n) = \frac{n!}{\theta_{(n)}} \prod_{i=1}^{n} \frac{\theta^{a_i}}{i^{a_i} a_i!},$$

where $\theta_{(n)} = \theta(\theta + 1) \cdots (\theta + n - 1)$. This formula has since played an important role in many different areas of probability theory, sometimes fairly distant from the original application to population genetics. Among many others, we refer the reader to [4] and to [28] for different proofs of (8).

Unfortunately, the methods used to prove (8) do not seem to apply to the more general framework of $\Lambda$-coalescents. In fact, there are very few explicit results studying the structure of a sample of the population in this context. Let us mention, in particular, the work of Möhle [39], Theorem 3.1, who gets a recursive formula for the allele frequency spectrum. However, this may be so intricate that this recursion is difficult to use in practice for moderately large sample sizes.

We present here an asymptotic formula for the frequency spectrum, both in the infinite alleles and the infinite sites models, as the sample size $n \to \infty$. We work under the convention that the genealogy of the population can be described by a $\Lambda$-coalescent $(\Pi_t, t \geq 0)$. We focus on the case where the measure $\Lambda$ is the $\text{Beta}(2 - \alpha, \alpha)$ distribution and $1 < \alpha < 2$. We assume that mutations occur at constant rate $\theta > 0$.

THEOREM 9. *Assume that $\Lambda$ has the $\text{Beta}(2-\alpha,\alpha)$ distribution with $1 < \alpha < 2$. Fix a positive integer $k$. Then*

$$n^{\alpha-2} M_k(n) \to_p \theta\alpha(\alpha-1)^2 \frac{\Gamma(k+\alpha-2)}{k!}$$

*and*

$$n^{\alpha-2} N_k(n) \to_p \theta\alpha(\alpha-1)^2 \frac{\Gamma(k+\alpha-2)}{k!},$$

*where $\to_p$ denotes convergence in probability as $n \to \infty$.*

REMARK 10. To understand where these results come from, recall that in Theorem 1.9 of [8], we showed that

$$(9) \qquad n^{\alpha-2} M(n) \to_p \theta \frac{\alpha(\alpha-1)\Gamma(\alpha)}{2-\alpha}.$$

In Section 5, we will show that for small times, the $\text{Beta}(2-\alpha,\alpha)$-coalescent can be approximately described by the genealogy of a continuous-time branching process in which individuals live for an exponential amount of time with mean 1 and then have a number of offspring distributed according to $\chi$, where $P(\chi=0) = P(\chi=1) = 0$ and where, for $k \geq 2$, we have

$$(10) \qquad P(\chi=k) = \frac{\alpha(2-\alpha)(3-\alpha)\cdots(k-1-\alpha)}{k!} = \frac{\alpha\Gamma(k-\alpha)}{k!\Gamma(2-\alpha)}.$$

This offspring distribution is supercritical with mean $1 + 1/(\alpha-1)$. We will show that if $\tau$ is an independent exponential random variable with mean $1/c$, where $c = (2-\alpha)/(\alpha-1) > 0$, and $k$ is a positive integer, then

$$(11) \qquad M_k(n) \sim M(n) P(\xi_\tau = k),$$

where

$$P(\xi_\tau = k) = \frac{(2-\alpha)\Gamma(k+\alpha-2)}{\Gamma(\alpha-1)k!}.$$

This result, and the analogous result for $N_k(n)$, will imply Theorem 9.

REMARK 11. One can only observe $M_k(n)$ from biological data if the ancestral type is known. Otherwise, it is necessary to work with the "wrapped frequency spectrum" $\hat{M}_k(n) = M_k(n) + M_{n-k}(n)$. For fixed $k \geq 1$, one can see from (11) that as $n \to \infty$, these two quantities have the same asymptotics because the limiting values of $M_k(n)/M(n)$ sum to one and therefore $M_{n-k}(n)/M(n)$ goes to zero in probability as $n \to \infty$.

REMARK 12. It is natural that the distribution (10) arises in this context because when the Beta$(2-\alpha,\alpha)$-coalescent has $b$ blocks, the probability that its next merger involves $k$ blocks converges to $P(\chi=k)$ as $b\to\infty$ (see [8, 13]). Of course, an individual having $k$ offspring in the Galton–Watson process corresponds to a merger of $k$ blocks in the corresponding coalescent process going backward in time.

REMARK 13. The limiting behavior described in Theorem 9 also arises in the theory of exchangeable partitions. Following Lemma 3.11 in Pitman [46], let $\Pi$ be an exchangeable random partition whose ranked asymptotic frequencies $P_i$ satisfy

$$(12) \qquad P_i \sim Zi^{-1/(2-\alpha)}$$

almost surely for some random variable $Z$ such that $0 < Z < \infty$. Then if $|\Pi_n|$ (resp., $|\Pi_{n,j}|$) is the number of blocks (resp., number of blocks of size $j$) of $\Pi$ restricted to $\{1,\ldots,n\}$, we have

$$(13) \qquad |\Pi_n| \sim S_\alpha n^{2-\alpha}$$

almost surely for a random variable $S_\alpha$ determined explicitly from $Z$. Moreover,

$$(14) \qquad \frac{|\Pi_{n,j}|}{|\Pi_n|} \sim \frac{(2-\alpha)\Gamma(k+\alpha-2)}{\Gamma(\alpha-1)k!}.$$

In fact, it follows from an unpublished work of Hansen and Pitman [26] that (13) implies (12), which, in turn, using Lemma 3.11 of [46], implies (14). Note that the distribution on the right-hand side of (14) previously appeared in the context of urn schemes in the work of Karlin [30]. See also [48], and see [43, 44], where this distribution occurs in the context of Brownian motion and related processes.

To connect these results to Theorem 9, let $\Pi$ be the allelic partition obtained by superimposing mutation marks at rate $\theta$ on the tree associated with a Beta-coalescent, started at time 0 with infinitely many individuals. Then $\Pi$ is an exchangeable partition and the restriction $\Pi_n$ of $\Pi$ to $\{1,\ldots,n\}$ has the same distribution as the partition described in Section 3.3. From this and (9), one can show (see, e.g., Lemma 34) that $n^{\alpha-2}|\Pi_n| \to_p s_\alpha$, where $s_\alpha$ is the constant on the right-hand side of (9). If one could show that this convergence holds almost surely, then this would supply an alternative proof of Theorem 9. Also, this would presumably work for coalescent processes satisfying the condition of Theorem 1.9 in [8]. However, we note that proving almost sure convergence is difficult due to the randomness of the asymptotic frequencies $P_i$.

## 4. The lookdown process in a continuous random tree.

4.1. *Branching processes obtained from superprocesses.* The lookdown process is a powerful tool introduced (and subsequently modified) by Donnelly and Kurtz [16] to encode the genealogy of a superprocess by a countable system of particles. We will describe it in a more general context than the one strictly needed for the applications we have in mind in this paper because we believe that this construction is of independent interest. However, the lookdown process can be defined even more generally than how we will do here (e.g., we will not treat the case where the particles are allowed to have some spatial motion and interact). The setting for this part is the following. We let $\psi$ be a branching mechanism with no Brownian component and no drift term, that is, there exists $a \in \mathbf{R}$ and a Lévy measure $\nu$ such that

$$(15) \qquad \psi(q) = aq + \int_0^\infty (e^{-qx} - 1 + qx\mathbf{1}_{\{x \leq 1\}})\nu(dx).$$

Rather than associating with $\psi$ a CSBP with this branching mechanism, we first construct a superprocess $M_t$ taking its values in the space of finite measures on $(0, 1)$, which is defined through its generator $L$: for a function $F$ acting on measures $\mu$ on $(0, 1)$,

$$(16) \quad LF(\mu) = a \int_0^1 \mu(dx) F'(\mu, x)$$
$$+ \int_0^1 \mu(dx) \int_0^\infty \nu(dh)(F(\mu + h\delta_x) - F(\mu) - \mathbf{1}_{\{h \leq 1\}} h F'(\mu, x)).$$

The notation $F'(\mu, x)$ stands for $\lim_{\varepsilon \to 0} \varepsilon^{-1}(F(\mu + \varepsilon \delta_x) - F(\mu))$ and accounts for an infinitesimal modification of $F$ in the direction $\delta_x$. If $\psi$ had a quadratic term, then there would be an extra term in the generator; see equation (1.15) in [14]. Note that for every $0 < r < 1$,

$$Z_t = M_t([0, r])$$

defines a $\psi$-CSBP started at $M_0([0, r])$. Indeed, applying the generator to a function $F(\mu) = \varphi(z)$, where $z = \mu([0, r])$, directly yields that the generator $L_1$ of the process $Z_t$ is

$$L_1 \varphi(z) = a \int_0^r \mu(dx) \varphi'(x)$$
$$+ \int_0^r \mu(dx) \int_0^\infty \nu(dh)(\varphi(z + h) - \varphi(z) - h\mathbf{1}_{\{h \leq 1\}} \varphi'(z))$$
$$= z L_2 \varphi(z)$$

since the second integral does not depend on $x$ and is equal to $L_2 \varphi(z)$, where $L_2$ is the generator of a Lévy process with Lévy exponent $\psi(q)$. By Lamperti's result relating a CSBP to a time-change of a Lévy process [35], we conclude that $Z_t$ is a $\psi$-CSBP. The interpretation of $M_t$ is as follows. If we imagine the population represented by $Z_t$ as a continuous population where each individual is endowed with an originally distinct label between 0 and 1 (and where individuals and their descendants have the same label), then $M_t([0,a])$ is the total number of individuals at time $t$ descending from some individual with a label between 0 and $a$. Another process of interest in this setting is the so-called *ratio process* $R_t = M_t/Z_t$, where $Z_t = M_t([0,1])$. Thus, for every $t$, $R_t$ is a probability distribution on $(0,1)$ which describes the composition of the population at a given time: the typical state at time $t > 0$ for $R_t$ (at least in the subcritical case, see below) is a linear combination of Dirac masses $\sum_i \rho_i \delta_{x_i}$, subject to $\sum_i \rho_i = 1$, where each atom corresponds to groups of individuals in the population at time $t$ descending from the same individual at time 0 (whose label was $x_i$) in proportion $\rho_i$.

4.2. *The lookdown process associated with a CSBP.* The purpose of the Donnelly–Kurtz construction is to give a representation of the ratio process $R_t$ as the limit of empirical distributions associated with a countable system of particles. A major consequence of this construction is a transparent notion of genealogy for $Z_t$, which is otherwise difficult to grasp in the context of a continuous population. What follows is largely inspired by [14] and [22], Chapter 5. To define the (modified) lookdown process, we have a countable number of individuals who will be identified with their *type*. Initially, individual $i$ has type $\xi_i(0)$. The types $\xi_i(0)$ for $i = 1, 2, \ldots$ are given by uniform i.i.d. random variables on $(0,1)$. At any given time $t$, $\xi_i(t)$ will be the type of the individual occupying level $i$. The variables $\xi_i(t)$ may change due to events called *birth events*. Suppose we have a countable configuration of space-time points,

$$n = \sum_i \delta_{(t_i, y_i)},$$

where $t_i \geq 0$ and $0 \leq y_i \leq 1$, and assume that $\sum_{t_i \leq t} y_i^2 < \infty$ for all $t \geq 0$. [Later, we will specify a point configuration $(t_i, y_i)$ associated with a CSBP.] Each atom $(t_i, y_i)$ corresponds to a birth event. At such a time, a proportion $y_i$ of levels is said to *participate* in the birth event: each level flips a coin with probability of heads $y_i$. Those which come up heads participate in the birth event. We describe the modification in the levels on the first $n$ levels. Suppose the levels participating are $1 \leq i_1 < i_2 < \cdots < i_k \leq n$. Then at time $t = t_i$, their type is modified by the following rule: for all $1 \leq j \leq k$, $\xi_{i_j}(t) = \xi_{i_1}(t^-)$. In other words, participating levels take the type of the smallest level participating. We do not destroy the individuals previously occupying levels $i_2, \ldots, i_k$, but, instead, we move $\xi_{i_2}(t^-)$ to the first level not participating in a birth event and keep shifting individuals upward, with each individual taking the first available spot. This is illustrated in Figure 4.

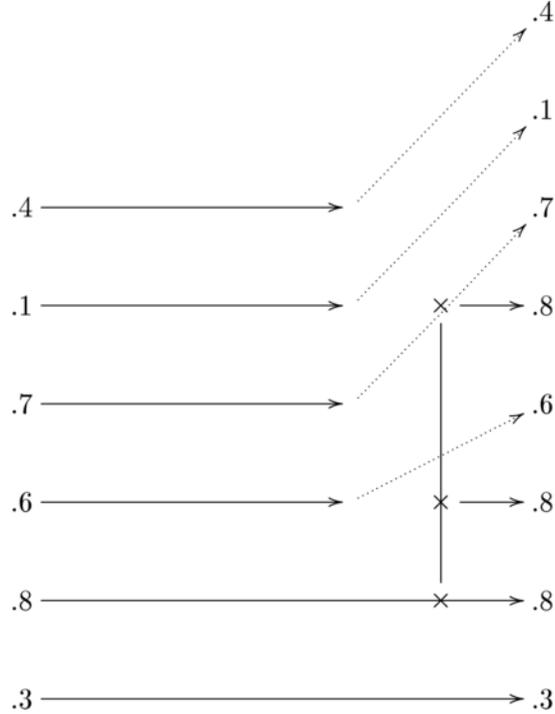

FIG. 4. *Representation of the lookdown process. Levels* 2, 4 *and* 5 *participate in a birth event. Other types get shifted upward. The numbers on the left and on the right indicate the types before and after the birth event, respectively.*

One way to make this construction rigorous is to observe that due to our assumption $\sum_{t_i \leq t} y_i^2 < \infty$, only finitely many birth events affect the first $n$ levels in any compact time-interval. The processes defined by this procedure are consistent by restriction as $n$ increases, so there is a well-defined process $(\xi_i(t), t \geq 0, i = 1, 2, \ldots)$ by Kolmogorov's extension theorem.

Having described the construction for a general configuration of space-time points $(t_i, y_i)$, we now restrict to the case where $(t_i, y_i)$ is given by the following construction. Let $Z_t(r)$ be a $\psi$-CSBP, where $\psi$ has the form (15) and where we have written the starting point $r > 0$ as an argument of $Z_t$. Let $\tau$ be the extinction time (which may not be finite a.s., but will be in the subcritical case in which we are interested). We only define the lookdown process until time $\tau^-$. With each time $t_i$ such that $\Delta Z_{t_i} > 0$, associate $y_i = \Delta Z_{t_i}/Z_{t_i}$ (observe that $0 \leq y_i \leq 1$). It is then standard to check that if $t < \tau$, then
$$\sum_{t_i \leq t} y_i^2 < \infty.$$
Indeed, one can bound $Z_{t_i}$ from below by $I_t = \inf_{0 \leq s \leq t} Z_s > 0$ so that this sum is smaller than
$$(I_t)^{-2} \sum_{t_i \leq t} (\Delta Z_{t_i})^2 < \infty$$
because $Z$ is obtained as a time-change of a Lévy process whose jumps are square-summable due to the fact that $\int_0^\infty (1 \wedge x^2)\nu(dx) < \infty$ and when $t < \tau$, the jumps of $Z$ are the jumps of the Lévy process in some random, but finite, time-interval.

Thus, there is a well-defined lookdown process $(\xi_i(t), t \geq 0, i = 1, 2, \ldots)$ associated with this sequence $(t_i, y_i)$. Observe that for all $t \geq 0$, $(\xi_i(t), i = 1, 2, \ldots)$ is an exchangeable sequence so that the limit
$$\rho_t = \lim_{n \to \infty} \frac{1}{n} \sum_{i=1}^\infty \delta_{\xi_i(t)}$$
is well defined by De Finetti's theorem. Then $(\rho_t, t \geq 0)$ has the same distribution as the process $(R_t, t \geq 0)$ obtained in the previous section from a superprocess $M_t$ started from $M_0 = r\mathbf{1}_{\{0 \leq x \leq 1\}} dx$ (see, e.g., the argument starting from (2.15) in [14]). To understand heuristically why this is true, note that when there is a jump in the CSBP, so $\Delta Z_t = x > 0$, some individual in the population has a large number of offspring, causing the proportion of individuals with the same type as this individual to have a jump of size $x/(Z_{t-} + x) = \Delta Z_t/Z_t$. This is precisely what happens in the lookdown process.

We now specialize to the subcritical case. That is, we assume that $\psi$ is a branching mechanism as in (15) and that
$$\int_1^\infty \frac{dq}{\psi(q)} < \infty.$$
By a well-known criterion of Grey [27], this ensures that $\tau < \infty$ a.s., that is, the population becomes extinct in finite time. Observe that one of the nontrivial features of the lookdown process is that since $Z_t$ becomes extinct in finite time, almost surely only finitely many individuals have descendants alive at time $t > 0$, which means that the composition of the population is made of finitely many different types of individuals and that, ultimately, only one type remains in the population. Note that this can happen in the lookdown process even though we never kill labels because some labels get pushed off to infinity due to the successive birth events, thus disappearing from the visible population. This feature will become apparent from our construction of the lookdown process in terms of the continuous random tree.

4.3. *Constructing the lookdown process from a continuous tree.* In this section, we will provide a construction of the lookdown process from a continuous random tree. Once again, we emphasize that the branching mechanism need not be stable. However, we will always assume subcriticality, that is,

$$\int_1^\infty \frac{dq}{\psi(q)} < \infty,$$

so that there is a continuous version of the height process and its local times are well defined (see [17]).

Before we start, we need to recall some facts about the height process. Associated with the process $H$ is an infinite measure $N$ which plays a role analogous to Itô's excursion measure for Brownian motion (see [47]). The excursion property for $(H_t, t \leq T_r)$ will be used on several occasions. It can be phrased as follows. Let $(g_i, d_i), i \in \mathcal{I}$ be the excursion intervals of $H$ above zero, so

$$\bigcup_{i \in \mathcal{I}} (g_i, d_i) = \{s \geq 0 : H_s > 0\}.$$

For each $i \in \mathcal{I}$, define the function $e_i$ by $e_i(s) = H_{g_i+s}$ for $0 \leq s \leq d_i - g_i$ and $e_i(s) = 0$ otherwise. Let $C_+([0,\infty))$ be the set of nonnegative real-valued functions defined on $[0,\infty)$. Recall that $(L_s^a, s \geq 0, a \geq 0)$ is the local time process for $H$. Then the random measure

$$\sum_{i \in \mathcal{I}} \delta_{(L_{g_i}^0, e_i)} \tag{17}$$

is a Poisson point process on $[0,\infty) \times C_+([0,\infty))$ with intensity measure $dl \times N(d\omega)$, where $dl$ denotes Lebesgue measure and $N(d\omega)$ is the excursion measure, which is a $\sigma$-finite measure on $C_+([0,\infty))$. More generally, $H$ (although not a Markov process in general) enjoys a similar excursion property above any given level $a > 0$. For each $a > 0$, let $(g_i^a, d_i^a), i \in \mathcal{I}^a$, be the connected components of the open set $\{s : H_s > a\}$. For each $i \in \mathcal{I}^a$, define the excursion $e_i^{(a)}$ by $e_i^{(a)}(s) = H_{g_i^a+s} - a$ for $0 \leq s \leq d_i^a - g_i^a$ and $e_i^{(a)}(s) = 0$ otherwise. For each $s \geq 0$, define

$$\widetilde{\tau}_s^a = \inf\left\{t : \int_0^t \mathbf{1}_{\{H_r \leq a\}} dr > s\right\}, \qquad \overline{\tau}_s^a = \inf\left\{t : \int_0^t \mathbf{1}_{\{H_r > a\}} dr > s\right\}.$$

Define the processes $(\widetilde{H}_s^a, s \geq 0)$ and $(\overline{H}_s^a, s \geq 0)$ such that $\widetilde{H}_s^a = H_{\widetilde{\tau}_s^a}$ and $\overline{H}_s^a = H_{\overline{\tau}_s^a} - a$. By Proposition 3.1 of [18], the random measure

$$\sum_{i \in \mathcal{I}^a} \delta_{(L_{g_i^a}^a, e_i^{(a)})} \tag{18}$$

is a Poisson point process on $[0,\infty) \times C_+([0,\infty))$ with intensity measure $dl \times N(d\omega)$ and is independent of $(\widetilde{H}_s^a, s \geq 0)$. Since $\overline{H}^a$ can be recovered from the random measure (18), a consequence of this result is that $(\overline{H}^a_s, s \geq 0)$ has the same law as $(H_s, s \geq 0)$ and is independent of $(\widetilde{H}^a_s, s \geq 0)$.

Having recalled this property, we now describe our construction of the lookdown process in a continuous random tree. Let $(Z_t, t \geq 0)$ be a $\psi$-CSBP started from $Z_0 = r > 0$ ($\psi$ is assumed to be subcritical) and assume that $Z_t$ is obtained as the local times of the height process $(H_t, t \leq T_r)$, as in (6). Let $\tilde{\xi} := ((\tilde{\xi}_j(t)), \ t \geq 0, j = 1, 2, \ldots)$ be a lookdown process obtained from $(Z_t, t \geq 0)$, as in the previous section. That is, it is obtained from the configuration of space-time points $(t_i, \Delta Z_{t_i}/Z_{t_i})$. The process $\tilde{\xi}$ will serve as a reference lookdown process to which we will compare the one we will construct below.

We will now construct a version $\xi$ of the process $\tilde{\xi}$ that will be entirely defined in terms of the height process $H$. We start by introducing some notation. Consider the height process $(H_t, t \leq T_r)$. The key point of this construction is that we choose a specific labeling for the excursions; namely, we rank the excursions according to their supremum. We denote by $e_j^{(t)}$ the $j$th highest excursion above the level $t$ (when $t = 0$, we sometimes simply write $e_j$ instead of $e_j^{(0)}$). We draw a sequence of i.i.d. random variables $(U_i)_{i \in \mathbb{N}}$ with the uniform distribution on $(0, 1)$. They will serve as the initial types in the lookdown construction, so that at any time, $\xi_j(t)$ is equal to one of the $U_i$'s. Thus, let $\xi_j(0) = U_j$ for all $j \geq 1$. Then for each $t > 0$, for each $j \geq 1$, we let $k(j, t)$ be the unique integer such that $e_j^{(t)}$, the $j$th highest excursion above $t$, is part of the excursion $e_{k(j,t)}^{(0)}$, the $k(j, t)$th highest excursion above 0, and we let

$$\xi_j(t) = U_{k(j,t)}.$$

We say that the excursion $e_j^{(t)}$ has type $U_{k(j,t)}$.

THEOREM 14. *The processes $\xi$ and $\tilde{\xi}$ have the same distribution. That is, $((\xi_j(t)), t \geq 0, j = 1, 2, \ldots)$ has the distribution of the modified lookdown construction associated with the CSBP $(Z_t, t \geq 0)$.*

Before we start proving this result, here is a description of the dynamics of the process $(\xi_j(t), t \geq 0)$. As $t$ increases, the relative ranking of the excursions above $t$ evolves. If $\Delta Z_t > 0$, then this means that with probability one, $H$ has (infinitely many) local minima at $t$, resulting in (infinitely many) additional excursions above $t$. Indeed, note that by Theorem 4.7 in [18], this corresponds to a unique excursion above $t^-$ splitting into infinitely many excursions. Moreover, all local minima of $(H_t, t \geq 0)$ are in fact associated with jumps of $Z_t$ (this would not be true if $\psi$ had a quadratic term; see Theorem 4.7 of [18]). We then say that some birth event happens. We rerank all excursions according to their new order (again, given by the rank of their supremum). Old excursions keep their old type (but might change their level) and the newly added excursions take the type from their father. If excursion $e_j^{(t)}$ splits, then this means that many levels $k$ with $k \geq j$ take the type $\xi_j(t)$. Those who do not take this type get shifted upward accordingly. To use the Donnelly–Kurtz terminology, we say that the levels $k \geq j$ adopting the type $\xi_j(t)$ *take part in the birth event*.

Let $\mathcal{F} = (\mathcal{F}_a, a \geq 0)$ be the filtration such that $\mathcal{F}_a = \sigma(\widetilde{H}^b, b \leq a)$. The key observation for the proof of Theorem 14 is summarized by the following lemma.

LEMMA 15. *Let $a > 0$ be a stopping time of the filtration $\mathcal{F}$ such that $\Delta Z_a > 0$ a.s. Define a sequence $(\varepsilon_i)_{i \in \mathbb{N}}$ by $\varepsilon_i = 1$ if the level $i$ takes part in the birth event at time $a$ for the process $\xi$ (i.e., the $i$th highest excursion above $a$ is a newly created excursion) and $0$ otherwise. Then the distribution of the sequence $(\varepsilon_i)_{i \in \mathbb{N}}$ is that of a sequence of i.i.d. Bernoulli variables with parameter $\Delta Z_a / Z_a$.*

PROOF. We know (see Theorem 4.7 in [18]) that if $\Delta Z_a > 0$, then $a$ is necessarily a level where exactly one excursion is splitting into infinitely many smaller ones (i.e., $a$ is a level where $H$ reaches a multiple infimum and for $b < a$, all those infima are reached within the same excursion above $b$). In other words, if $a$ is a jump time of $Z$, there is a unique interval $(s, t)$ such that $L_t^{a-} = L_s^{a-}$ and $L_t^a - L_s^a = \Delta Z_a$. Let us denote $x = L_s^a$ and $y = L_t^a$.

For $i \geq 1$, define $h_i^{(a)} := \max e_i^{(a)}$ to be the height of the $i$th highest excursion above level $a$ and let $t_i^{(a)}$ denote the local time accumulated at level $a$ when the excursion $e_i^{(a)}$ starts. By applying the strong Markov property which will be proved at the very end of this section, in Lemma 17 we see that conditionally on $Z_a$, the process $\overline{H}_t^a$ has the same distribution as $H$ run until $T_{Z_a}$. Hence, the atoms $(t_i^{(a)}, h_i^{(a)})$ form a Poisson point process on $[0, Z_a] \times \mathbb{R}^+$ with intensity measure $dt \times n(dh)$, where $n$ is absolutely continuous with respect to the Lebesgue measure, $n(0, \infty) = \infty$ and $n(h, \infty) < \infty$ for $h > 0$. The measure $n$ is the "law" of the heights of excursions under the measure $N$.

Observe that the levels that take part in the birth event are exactly the levels $k$ which correspond to the rank of a newly created excursion $e_k^{(a)}$, that is, the excursion such that $t_k^{(a)} \in (x, y)$, where $(x, y)$ is the new interval of local time. The statement then amounts to the well-known fact about Poisson point processes that the $t_j$ (observe that $t_j$ is the time of the $j$th record of the Poisson point process) are i.i.d. uniformly distributed random variables over $(0, Z_a)$ and are independent of the sequence of records $h_j^{(a)}$. As the events $\{t_j \in (x, y)\}$ and $\{\varepsilon_j = 1\}$ coincide, the conclusion follows. $\square$

Now, fix $\varepsilon > 0$. Let $a_1$ be the first time $t$ such that $\Delta Z_t / Z_t > \varepsilon$. Observe that almost surely $a_1 > 0$ and that $a_1$ is a stopping time for $\mathcal{F}$. We may thus define, inductively, $a_1 < a_2 < \ldots$, the set of stopping times such that $\Delta Z_t / Z_t > \varepsilon$ and for each $i \geq 1$, $a_i$ is a stopping time of $\mathcal{F}$. For $i \geq 1$, a multiple infimum is reached at level $a_i$, which corresponds to a single excursion that splits into an infinite number of descendants at this precise level. Define a process $(\xi_j^{(\varepsilon)}(t), t \geq 0, j = 1, 2, \ldots)$ as follows:

- if $t$ is not a jump time for $Z$, then nothing happens for $\xi^{(\varepsilon)}$, that is, we have $\xi^{(\varepsilon)}(t^-) = \xi^{(\varepsilon)}(t)$;
- if $t$ is a jump time for $Z$, but $\Delta Z_t / Z_t < \varepsilon$, we use an independent coin flipping with probability of heads $y = \Delta Z_t / Z_t$, and the standard Donnelly–Kurtz procedure, to obtain $\xi^{(\varepsilon)}(t)$ from $\xi^{(\varepsilon)}(t-)$;
- if $t$ is a jump time for $Z$ and $\Delta Z_t / Z_t \geq \varepsilon$ (i.e., $t = a_i$ for some $i$), we say that the levels which take part in the birth event are exactly the relative ranks of the newly created excursions at level $t$.

LEMMA 16. *For each fixed $\varepsilon > 0$, the processes $\xi^{(\varepsilon)}$ and $\tilde{\xi}$ have the same distribution.*

PROOF. We only need to show that our new rule for the times $a_i$ does not differ from the usual construction. As the $a_i$'s are a sequence of stopping times, we can apply Lemma 15 to see that we are again deciding who takes part in the birth event according to a sequence of i.i.d. Bernoulli variables with parameters $\Delta Z_{a_i} / Z_{a_i}$. The strong Markov property also implies that the sequences used at the successive times $a_i$ are independent. Hence, $\xi^{(\varepsilon)}$ has the same distribution as $\tilde{\xi}$. $\square$

PROOF OF THEOREM 14. Let $b_1, \ldots, b_m$ be the times at which there is a change in the first $n$ levels for the process $\xi$ (the number $m$ of such times is necessarily at most $n - 1$ since at each of the $b_i$, the diversity of types among the first $n$ levels must be reduced at least by 1). Let $F$ be a bounded functional on the Skorokhod space $\mathbb{D}(\mathbb{R}_+^\infty, \mathbb{R})$ endowed with the product topology inherited from $\mathbb{D}(\mathbb{R}_+, \mathbb{R})$ and assume that $F$ only depends on the first $n$ coordinates (levels) for some arbitrarily fixed number $n \geq 1$. Then

$$(19) \quad |E(F(\xi)) - E(F(\xi^{(\varepsilon)}))| \leq \|F\|_\infty P(\{b_1, \ldots, b_m\} \not\subset \{a_1, a_2, \ldots\})$$

because when $\{b_1, \ldots, b_m\} \subset \{a_1, a_2, \ldots\}$, the first $n$ coordinates of $\xi^{(\varepsilon)}$ and $\xi$ coincide exactly. Since $\xi^{(\varepsilon)}$ and $\tilde{\xi}$ have the same distribution, by Lemma

16, we deduce that

(20) $\quad |E(F(\xi)) - E(F(\tilde{\xi}))| \leq \|F\|_\infty P(\{b_1, \ldots, b_m\} \not\subset \{a_1, a_2, \ldots\}).$

Note that

$$\lim_{\varepsilon \to 0} P(\{b_1, \ldots, b_m\} \not\subset \{a_1, a_2, \ldots\}) = 0.$$

Indeed, there are only finitely many jumps affecting the first $n$ levels, so

$$\eta := \inf_{t \in \{b_1, \ldots, b_m\}} \frac{\Delta Z_t}{Z_t} > 0 \qquad \text{a.s.}$$

Since $\{b_1, \ldots, b_m\} \not\subset \{a_1, a_2, \ldots\}$ is equivalent to $\eta < \varepsilon$, we see that

$$P(\{b_1, \ldots, b_m\} \not\subset \{a_1, a_2, \ldots\}) = P(\eta < \varepsilon) \to 0$$

as $\varepsilon \to 0$ because $\eta > 0$ a.s. It follows by letting $\varepsilon \to 0$ in (20) that the restrictions of $\xi$ and $\tilde{\xi}$ to the first $n$ coordinates are identical in distribution. By the uniqueness in Kolmogorov's extension theorem, the processes $\xi$ and $\tilde{\xi}$ are thus identical in distribution. $\square$

It now remains to establish the strong Markov property, which we used on several occasions. Note that this lemma holds even at stopping times $T$ such that $\Delta Z_T > 0$.

LEMMA 17. *Let $T$ be a stopping time of $\mathcal{F}$. Conditionally on $Z_T = z$, the processes $\overline{H}_t^T$ and $\widetilde{H}_t^T$ are independent. Moreover, $\overline{H}_t^T$ is distributed as $(H_t, t \leq T_z)$.*

PROOF. When $T = s$ is a deterministic stopping time, then this is the content of Corollary 3.2 in [18]. Suppose we now try to verify the claim when $T$ is a stopping time of $Z$ which can only take a countable number of values $\{t_k\}$, say. Let $F$, $G$ be two nonnegative functions defined on $C([0, \infty])$ and assume that they are continuous for the topology of uniform convergence on compact sets. Since $\{T = t_k\}$ is $\mathcal{F}_{t_k}$-measurable we then have

$$\begin{aligned}
&E[F(\widetilde{H}_t^T, t \geq 0) G(\overline{H}_t^T, t \geq 0) | Z_T = z] \\
&\quad = \sum_{k \geq 0} E[F(\widetilde{H}_t^{t_k}, t \geq 0) G(\overline{H}_t^{t_k}, t \geq 0) \mathbf{1}_{\{T=t_k\}} | Z_{t_k} = z] \\
&\quad = \sum_{k \geq 0} E[G(H_{t \wedge T_z}, t \geq 0)] E[F(\widetilde{H}_t^T, t \geq 0) \mathbf{1}_{\{T=t_k\}} | Z_T = z] \\
&\quad = E[G(H_{t \wedge T_z}, t \geq 0)] E[F(\widetilde{H}_t^T, t \geq 0) | Z_T = z].
\end{aligned}$$

To extend this to stopping times taking a continuous set of values, we use standard approximations of the stopping time $T$ by

$$T_n = \sum_{k \geq 0} \frac{k+1}{2^n} \mathbf{1}_{\{k/2^n \leq T < (k+1)/2^n\}}.$$

Note that $T_n$ approaches $T$ from above within $2^{-n}$. To begin observe that

$$\int_0^{T_r} \mathbf{1}_{\{T \leq H_u \leq T_n\}} \, du = \int_T^{T_n} Z_a \, da,$$

which, by (right) continuity of $Z$ at $T$, is smaller than $C2^{-n}$ for $n$ sufficiently large a.s. To see that $\overline{H}_s^{T_n}$ approaches uniformly $\overline{H}_s^{T}$, we think of the following picture. There are two sources of difference between $\overline{H}_s^{T_n}$ and $\overline{H}_s^{T}$. One is a shift downward for the excursions above 0 because the parts of an excursion between $T$ and $T_n$ are erased in $\overline{H}_t^{T_n}$. This shift is at most $2^{-n}$. The other source is that there may be some excursions above $T$ that are not counted as excursions above $T_n$, or an excursion above $T$ could be split into two or more excursions above $T_n$ because of a local minimum between $T$ and $T_n$. This results in a horizontal shift. The total duration of this horizontal shift may never exceed the total time spent by $H$ in the strip $[T, T_n]$, which is not more than $C2^{-n}$, by the above remark. Hence, by uniform continuity of $H$, $\overline{H}_s^{T_n}$ approaches uniformly $\overline{H}_s^{T}$. A moment's thought shows that the same reasoning applies to $\widetilde{H}_s^{T_n}$ (and this does not require left continuity of $Z$ at $T$).

Therefore, if $F, G$ are, as above, two bounded, nonnegative and continuous functions on $C([0, \infty])$ and if $\varphi$ is also a bounded, continuous, nonnegative function on $\mathbb{R}$, since $T_n$ is a stopping time that takes only countably many values, we have

$$E[F(\widetilde{H}_t^{T_n}, t \geq 0) G(\overline{H}_t^{T_n}, t \geq 0) \varphi(Z_{T_n})]$$
$$= \int_0^\infty P(Z_{T_n} \in dz) \varphi(z) E[G(H_{t \wedge T_z}, t \geq 0)] E[F(\widetilde{H}_t^{T_n}, t \geq 0) | Z_{T_n} = z].$$

If $H'$ is another height process, independent of everything else, and if $L_n = \inf\{t > 0, L_t^0(H') > Z_{T_n}\}$, this can be rewritten as

$$E[F(\widetilde{H}_t^{T_n}, t \geq 0) G(\overline{H}_t^{T_n}, t \geq 0) \varphi(Z_{T_n})]$$
$$= E[F(\widetilde{H}_t^{T_n}, t \geq 0) G(H'_{t \wedge L_n}, t \geq 0) \varphi(Z_{T_n})].$$

Note that if $L = \inf\{t > 0, L_t^0(H') > Z_T\}$, then $(H'_{t \wedge L_n}, t \geq 0) \to (H'_{t \wedge L}, t \geq 0)$ uniformly almost surely. Indeed, because $Z_T$ is independent of $H'$, it suffices to show, by Fubini's theorem, that for a given $z$, $(H'_{t \wedge T'_{z \pm \varepsilon}}, t \geq 0)$ converges uniformly almost surely as $\varepsilon \to 0$ to $(H'_{t \wedge T'_z}, t \geq 0)$, where $T'$ is the inverse local time at 0 of $H'$. To see this, dropping the prime from the notation, first note that $T.$ is continuous at $z$ almost surely because it is a subordinator and, as such, does not have fixed discontinuities. Moreover, note that $\sup_{T_z \leq s \leq T_{z+\varepsilon}} H_s$, say, is the supremum of the heights of the excursions between $T_z$ and $T_{z+\varepsilon}$. By the excursion theory for $H$, this can be written as $S_\varepsilon = \sup_{t_i \leq \varepsilon} h(e_i)$, where $(t_i, h(e_i))$ is the Poisson point process of the heights of the excursions on an interval of duration $\varepsilon$. For any $\delta > 0$, excursions of height greater than $\delta$ have finite measure under $N$ and therefore $S_\varepsilon \leq \delta$ for sufficiently small $\varepsilon$. It follows that $S_\varepsilon \to 0$ as $\varepsilon \to 0$ almost surely or, in other words,

$$\|H_{t \wedge T_{z+\varepsilon}} - H_{t \wedge T_z}\|_\infty \to 0$$

almost surely. Therefore, $(H_{t \wedge T_{z \pm \varepsilon}}, t \geq 0)$ converges uniformly to $(H_{t \wedge T_z}, t \geq 0)$ a.s.

Since, on the other hand, $\widetilde{H}_t^{T_n}$ converges uniformly to $\widetilde{H}_t^T$ a.s., and since, similarly, $\overline{H}^{T_n}$ converges a.s. uniformly to $\overline{H}^T$ in the left-hand side, we conclude, by Lebesgue's dominated convergence theorem, that

$$E[F(\widetilde{H}_t^T, t \geq 0) G(\overline{H}_t^T, t \geq 0) \varphi(Z_T)]$$
$$= E[F(\widetilde{H}_t^T, t \geq 0) G(H'_{t \wedge T'_{Z_T}}, t \geq 0) \varphi(Z_T)].$$

From this, we immediately deduce, by conditioning on $Z_T = z$, the desired identity,

$$E[F(\widetilde{H}_t^T, t \geq 0) G(\overline{H}_t^T, t \geq 0) | Z_T = z]$$
$$= E[G(H_{t \wedge T_z}, t \geq 0)] E[F(\widetilde{H}_t^T, t \geq 0) | Z_T = z]. \qquad \square$$

4.4. *Proofs of Theorem* 1 *and Corollary* 2.

PROOF OF THEOREM 1. By Theorem 2.1 in [14], the time-changed genealogy of $Z_{R^{-1}(t)}$, as defined from the lookdown process, is a Beta$(2-\alpha, \alpha)$-coalescent. It then suffices to show that the notion of genealogy as we have defined it from the height process coincides with the notion of genealogy for the lookdown process constructed on the CRT.

There is a natural notion of genealogy associated with the lookdown construction. Namely, for any pair $i, j \geq 1$ and any times $0 \leq t \leq T$, we can decide if the levels $i$ and $j$ at time $T$ descend from the same level at time $t$ (more precisely, we can track their labels by going backward from time $T$ to time $t$ to see if they come from the same label).

When the lookdown construction is obtained (as explained above) from the process $H$, this means that levels $i$ and $j$ at time $T$ have the same ancestor at time $t$ if and only if the $i$th and $j$th highest excursions above $T$ are descendants of the same excursion above $t$.

Recall that $(V_i, i = 1, 2, \ldots)$ is a sequence of variables in $[0, T_r]$ where each $V_i$ is the left endpoint of the $i$th highest excursion above $R^{-1}(t)$. It is clear that if two excursions $e_i^{(R^{-1}(t))}$ and $e_j^{(R^{-1}(t))}$ above $R^{-1}(t)$ descend from the same excursion above $s$, then $V_i$ and $V_j$ are straddled by this excursion above $s$ or, in other words, that $\min_{r \in (V_i, V_j)} H(r) > s$. Hence, we see that the partition-valued process $(\Pi(s), 0 \le s \le t)$ such that $i$ and $j$ are in the same block of $\Pi(s)$ if and only if $\min_{r \in (V_i, V_j)} H(r) > R^{-1}(t-s)$, is exactly the process of the ancestral partition of the lookdown process $\xi$ between times $R^{-1}(t)$ and $R^{-1}(t-s)$. By applying Theorem 2.1 in [14], this entails that when $H$ is the height process associated with the $\alpha$-stable branching mechanism, $\Pi$ is a Beta$(2-\alpha, \alpha)$-coalescent—this was the content of Theorem 1. □

PROOF OF COROLLARY 2. Again, observe that the genealogy as defined from the lookdown process coincides with the following definition: $i$ and $j$ are in the same block of $\Pi_s$ if the $i$th and the $j$th highest excursions above level $R^{-1}(t)$ are subexcursions of a single excursion above $R^{-1}(t-s)$. Let $N_s$ be the number of excursions between $R^{-1}(t-s)$ and $R^{-1}(t)$ and, conditionally on $N_s = k$, number these excursions in random order $e_1, \ldots, e_k$, and let $\ell_1, \ell_2, \ldots, \ell_k$ be their respective local times at $R^{-1}(t)$. We want to show that the asymptotic frequency of the block corresponding to an excursion is proportional to $\ell$. However, reasoning as in Lemma 15, we see that, conditionally on $N_s = k$ and conditionally on $\ell_1, \ell_2, \ldots, \ell_k$, each level $i$ in the lookdown process at time $R^{-1}(t)$ falls in excursion $i$ with a probability that is equal to $\ell_i / Z_{R^{-1}(t)}$. It follows immediately from the law of large numbers that the asymptotic frequency of the block associated with $e_i$ is $\ell_i / Z_{R^{-1}(t)}$. In other words, the sequence of ranked frequencies of the ancestral partition defined by the lookdown process is almost surely equal to the process $(X(s), 0 \le s \le t)$. Corollary 2 immediately follows. □

**5. Small-time behavior and multifractal spectrum.** In this section, we use Theorem 1 to prove Theorems 4 and 5. We start by introducing our main tool, reduced trees.

5.1. *Reduced trees as Galton–Watson processes.* The key ingredient for the theorems in this section is the *reduced tree* associated with a height process $H$. For a fixed level $a$, the *reduced tree at level $a$* is a tree such that the number of branches of the tree at height $0 \le t \le 1$ is the number of excursions of $H$ above level $at$ that reach level $a$, with the natural genealogical structure defined by saying that $v$ is an ancestor of $w$ if the excursion associated with $w$ is contained in $v$. We will deduce from results of Duquesne and Le Gall [17] that when $H$ is the height process associated with the $\alpha$-stable branching mechanism, this tree is a Galton–Watson tree whose reproduction law can be described explicitly.

When the Beta-coalescent is constructed from the continuous random tree, the number of blocks $N(s)$ at time $s$ corresponds to the number of excursions above level $s'$ that reach level $R^{-1}(t)$, for some $s'$ and $t$. We can deduce the limiting behavior of $N(s)$ when $s \to 0$ from the limiting behavior of the reduced tree as $s' \to R^{-1}(t)$. However, because the reduced tree is a Galton–Watson tree, its limiting behavior is described by the Kesten–Stigum theorem, as stated in (22) below, and this leads to a proof of Theorem 4. Likewise, Theorem 5 is established by relating the multifractal spectrum of Beta-coalescents to the multifractal spectrum for Galton–Watson trees and then applying recent results of Mörters and Shieh [41] on the branching measure of Galton–Watson trees. An important step in the proof of these theorems is showing that events concerning the reduced tree at a fixed level can be carried over to the reduced tree at the random level $R^{-1}(t)$.

We now introduce more carefully the concept of reduced trees. We start with some notation. If $u > 0$, let $N_{(u)}$ denote the excursion measure of the height process, conditioned to hit level $u$,

$$N_{(u)}(\cdot) = N\left(\cdot \Big| \sup_{s \geq 0} H_s > u\right),$$

which is well defined since $N(\sup_{s \geq 0} H_s > u) < \infty$ for all $u > 0$. Let $(H_s, s \leq \zeta)$ be a realization of $N_{(u)}$ and consider the process $(\theta^u(t), 0 \leq t \leq u)$ defined by $\theta^u(t) = \# \mathrm{exc}_{t,u}$, the number of excursions above level $t$ reaching $u$ of $H$. Simple arguments show that almost surely for all $t < u$, we have $\theta^u(t) < \infty$.

DEFINITION 18. The reduced tree $\mathbf{T}^u$ at level $u$ associated with $(H_s, s \leq \zeta)$ is the tree encoded by the process $(\theta^u(tu), 0 \leq t \leq 1)$. In other words, each branch at level $0 \leq t \leq 1$ is associated with a unique excursion above level $tu$ reaching $u$.

In the context of quadratic branching where the height process is reflecting Brownian motion, this is a variant of a process already considered by Neveu and Pitman [42]. We should emphasize that, by a slight abuse of notation we will sometimes use the notation $\mathbf{T}^u$ even when the underlying process $(H_s, s \leq \zeta)$ is not a realization of $N_{(u)}$, but, rather the height process considered until time $T_r$, where it has accumulated local time $r$ at zero. In this case, $\mathbf{T}^u$ is, in fact, a forest consisting of a Poisson number of independent realizations of the tree of Definition 18. The following fact will be a crucial tool for much of our analysis. It states that up to a deterministic exponential time-change, the tree $\mathbf{T}^u$ is a continuous-time supercritical Galton–Watson (discrete) tree. We recall that here the branching mechanism is assumed to be stable.

PROPOSITION 19. *For fixed $u > 0$, the process $(\theta^u(u(1 - e^{-t})), 0 \leq t < \infty)$ is a continuous-time Galton–Watson process where individuals reproduce at rate 1 with a number of offspring $\chi$ satisfying*

$$E(r^\chi) = \frac{(1-r)^\alpha - 1 + \alpha r}{\alpha - 1}. \tag{21}$$

*More explicitly,*

$$P(\chi = k) = \frac{\alpha(2-\alpha)(3-\alpha)\cdots(k-1-\alpha)}{k!}, \qquad k \geq 2,$$

*and $P(\chi = k) = 0$ for $k \in \{0, 1\}$.*

PROOF. We show how this result follows from a result in Duquesne and Le Gall [17]. To simplify, we will assume that $u = 1$. By the remark following Theorem 2.7.1 of [17], the time of the first split $\gamma$ in $\theta^1(t)$ is a uniform random variable on $(0, 1)$. Then, conditionally on $\gamma = t$ and $\theta^1(\gamma) = k$, the process $Z_{\gamma+s}$ is distributed as the sum of $k$ independent copies of $(\theta^{1-t}(s), 0 \leq s \leq 1 - t)$. In particular, if we follow a branch in the tree from level 0 to level 1, we see that the times at which the corresponding individual reproduces are distributed according to the standard "stick-breaking" construction of a Poisson–Dirichlet random variable, described as follows. A first cut point is selected uniformly at random in $(0, 1)$ and the left piece is discarded. Another point is selected uniformly in the right piece. Discarding the left piece, we proceed further by selecting a point uniformly in the piece left after the second cut, and so on. It is well known and easy to see that the image of these points by the map $t \mapsto -\ln(1 - t)$ is a standard Poisson process. The distribution of the number of offspring at each branch point is naturally given by the law of the random variable $\theta^1(\gamma)$, whose distribution is identified in the remark following Theorem 2.7.1 of [17]. This implies the proposition. $\square$

REMARK 20. We also present an intuitive, but less precise, argument for why $(\theta^1(1 - e^{-t}); t \geq 0)$ is a Galton–Watson process (in the case of a stable branching mechanism). We recall that the process $\overline{H}_t^a$ is independent of the process $\widetilde{H}^a$, conditionally given the local time at level $a = 1 - e^{-t}$, and the excursions are given by the points of a Poisson point process with intensity $dl \times N(de)$; see (18) for a precise formulation. In particular, given that $k$ of them reach level 1, they are $k$ independent realizations of $N_{(e^{-t})}$. This proves the independence of $\theta^1(1 - e^{-(t+s)})$ with respect to its past, conditionally given $\theta^1(1 - e^{-t})$. Moreover, the law of each of these $k$ subtrees is identical to that of the whole tree. Indeed, the descendants at level $1 - e^{-(t+s)}$ of some excursion above $1 - e^{-t}$ reaching 1 is identical in law, after scaling the vertical axis by $e^{-t}$, to the descendants at level $1 - e^{-s}$ of an excursion above level 0 reaching level 1. [Recall that because the branching mechanism is stable, the height process has the following scaling property: if $(H_s, s \geq 0)$ is the height process under the measure $N_{(1)}$, then $H^{(u)} = (uH_{su^{-\alpha}}, s \geq 0)$ is a realization of $N_{(u)}$.] This proves that $|\mathbf{T}^1(1 - e^{-t})|$ is a Galton–Watson process. Observe, however, that this scaling argument does not give the reproduction rate of individuals, nor the exact offspring distribution.

We conclude this section by observing that the Galton–Watson process $[\theta^u(u(1 - e^{-t})), t \geq 0]$ satisfies the conditions needed to apply the celebrated Kesten–Stigum theorem. More precisely, we have the following lemma.

LEMMA 21. *There exists a random variable $W$ with $W > 0$ almost surely such that*

$$(22) \qquad e^{-t/(\alpha-1)}\theta^u(u(1-e^{-t})) \to W \qquad \text{a.s. when } t \to \infty.$$

PROOF. It can be checked that the reproduction law $\chi$ has mean $m = 1 + 1/(\alpha - 1)$. The Galton–Watson process is thus supercritical. Moreover, $P(\chi \geq k)$ decays like $k^{-\alpha}$ and in particular, $E(\chi \log \chi) < \infty$, so we may apply to this supercritical Galton–Watson process the Kesten–Stigum theorem (in continuous time) [5], Theorem 7.1. □

5.2. *Proof of Theorem 4 (number of blocks)*. We will now show that the variable $W$ in (22) above is a quantity which can be expressed in terms of the local time at level $u$. We start by focusing on the case $u = 1$ and we work under the measure $N_{(1)}$.

First, we need a simple continuity lemma for the local time at level 1 under $N_{(1)}$. Let $Z_t^{(1)}$ denote the total local time of the process $H$ at level $t$.

LEMMA 22. *Under $N_{(1)}$, $Z_t^{(1)}$ is continuous at $t = 1$, that is, $Z_{1-}^{(1)} = Z_1^{(1)}$, $N_{(1)}$-a.s.*

PROOF. When $Z_t$ is the local time at level $t$ of $(H_s, s \leq T_r)$, then it is well known that $Z_t$ cannot have a discontinuity at level 1 (indeed, $Z_t$ is a CSBP started at $Z_0 = 1$, hence it is a Feller process and so cannot have a fixed discontinuity). Conditionally on

$$\# \text{exc}_{0,1} = 1,$$

the excursion that reaches 1 is a realization of $N_{(1)}$ and as $\# \text{exc}_{0,1}$ is Poissonian, this event has strictly positive probability. Hence, the result follows. □

We now give the interpretation of $W$ in terms of $Z_1^{(1)}$.

LEMMA 23. *Let $K = (\alpha - 1)^{-1/(\alpha-1)}$ and let $u > 0$. Under $N_{(u)}$, we have*

$$\varepsilon^{1/(\alpha-1)} \theta^u(u(1-\varepsilon)) \to K u^{-1/(\alpha-1)} Z_u^{(u)} \qquad a.s.$$

*as $\varepsilon \to 0$, where $Z_u^{(u)}$ denotes the local time of $H$ at level $u$.*

REMARK 24. This result is thus a generalization of Lévy's result for the local time of Brownian motion as the limit of the rescaled "downcrossing number" (see, e.g., [47]). A similar result on the upcrossing number also exists and is, in fact, much simpler than the one we prove here due to the existence of an excursion theory *above* a fixed level.

PROOF. For simplicity, we will prove this result assuming that $u = 1$, but the case of general $u$ follows exactly the same arguments. We thus wish to prove that

$$\varepsilon^{1/(\alpha-1)} \theta^1(1-\varepsilon) \to K Z_1^{(1)} \qquad \text{a.s.},$$

as $\varepsilon \to 0$. We already know, by Lemma 21, that $\varepsilon^{1/(\alpha-1)} \theta^1(1-\varepsilon)$ converges almost surely to $W$. Hence, it is enough to prove the convergence in probability here to obtain that $W = K Z_1^{(1)}$ a.s. and to thereby conclude.

By excursion theory, conditionally on $Z_{1-\varepsilon}^{(1)} = z_\varepsilon$, the number of excursions above $1-\varepsilon$ that reach 1 is Poisson distributed with mean $z_\varepsilon N(\sup_{s \geq 0} H_s > \varepsilon)$. Now, recall that by [17], Corollary 1.4.2 applied with $\psi(u) = u^\alpha$,

$$N\left(\sup_{s \geq 0} H_s > \varepsilon\right) = (\alpha - 1)^{-1/(\alpha-1)} \varepsilon^{-1/(\alpha-1)} = K \varepsilon^{-1/(\alpha-1)}$$

[this is why the factor $u^{-1/(\alpha-1)}$ appears in the limit when $u \neq 1$ since, in this case, we need to compute $N(\sup_{s \geq 0} H_s > u\varepsilon)$]. Let $\delta > 0$ and let us show that

(23) $$P(\varepsilon^{1/(\alpha-1)} \theta^1(1-\varepsilon) > K Z_1^{(1)}(1+\delta)) \to 0$$

as $\varepsilon \to 0$. To do this, note that this is smaller than

(24) $$\begin{aligned} & P(|Z_1^{(1)} - Z_{1-\varepsilon}^{(1)}| > K Z_1^{(1)} \delta/2) \\ & \quad + P(\varepsilon^{1/(\alpha-1)} \theta^1(1-\varepsilon) > K Z_{1-\varepsilon}^{(1)}(1+\delta/2)). \end{aligned}$$

The first term converges to 0 by continuity of $Z$ at level 1. On the other hand, Markov's inequality implies that if $X$ is a Poisson random variable with mean $m/\varepsilon$, then for every $\lambda > 0$,

$$P(\varepsilon X > m(1+x)) \leq \exp\left[\frac{m}{\varepsilon}(-1 + e^\lambda - \lambda(1+x))\right].$$

By choosing $\lambda > 0$ sufficiently close to 0, we can find $c > 0$ such that

$$P(\varepsilon X > m(1+x)) \leq \exp(-cm/\varepsilon).$$

Therefore, the second term in (24) is bounded from above by

$$E(\exp(-c'Z^{(1)}_{1-\varepsilon}\varepsilon^{-1/(\alpha-1)})) \to 0$$

for some $c' > 0$, by Lebesgue's dominated convergence theorem, since $Z^{(1)}_{1-\varepsilon} \to Z^{(1)}_1$, a.s. This gives the convergence in probability for the lemma. $\square$

We note that the case $u \neq 1$ can also be obtained from the case $u = 1$ by using scaling properties of the process $H$: if $(H_s, s \geq 0)$ is the height process under the measure $N_{(1)}$, then $H^{(u)} = (uH_{su^{-\alpha/(\alpha-1)}}, s \geq 0)$ is a realization of $N_{(u)}$ (see, e.g., the remark before Theorem 3.3.3 of [17]).

LEMMA 25. *Assume that $\theta^u(t)$ is obtained for $0 \leq t \leq u$ from the reduced tree associated with the process $(H_t, 0 \leq t \leq T_r)$. Then*

$$(25) \qquad \lim_{t \to 0} t^{1/(\alpha-1)} \theta^u(1-t) \to K u^{-1/(\alpha-1)} Z_u \qquad a.s.$$

PROOF. This is a simple extension of Lemma 23. Again, to simplify, assume that $u = 1$. There is a slight difference with Lemma 23, because this was stated under the measure $N_{(1)}$, whereas here, $\mathbf{T}^1$ is defined from the height process $(H_s, s \leq T_r)$ and not a realization of $N_{(1)}$. However, this does not change the limit result, since the excursions of $(H_s, s \leq T_r)$ reaching 1 are independent and distributed with law $N_{(1)}$ (note that the result is trivially true when no excursion reaches level 1). Therefore, the result remains the same. $\square$

The point of the next lemma is to show that any almost sure property $A_u$ of the tree $\mathbf{T}^u$ still holds almost surely when the fixed level $u$ is replaced by the random level $R^{-1}(t)$ if we choose $t$ outside a deterministic set of Lebesgue measure 0. By convention, if $\mathbf{T}^u$ is empty (i.e., if $\sup_{0 \leq s \leq T_r} H_s < u$), we declare any property to be true by default. Since we wish to study the property $A_u$ at level $u = R^{-1}(t)$ for some $t$ and $\mathbf{T}^u$ is never empty, this will never play any role.

LEMMA 26. *Let $A_u$ be a property of the tree $\mathbf{T}^u$ such that for every $u > 0$, $P(A_u | \sup_{0 \leq s \leq T_r} H_s > u) = 1$. Then the set of $t$ such that $P(A_{R^{-1}(t)}) < 1$ has zero Lebesgue measure.*

PROOF. Let $F$ be the set of $t$ such that $A_t$ fails. By Fubini's theorem,

$$E \int_0^\infty \mathbf{1}_{\{t \in F\}} \, dt = 0.$$

Therefore, $\text{Leb}(F) = 0$ a.s. On the other hand, $t \mapsto R(t)$ is almost surely an absolutely continuous function. Indeed, it has a derivative at all points where $Z$ is continuous and $Z$ has only countably many discontinuities a.s. Therefore, $R(F)$ also has zero Lebesgue measure almost surely. Hence,

$$\int_0^\infty \mathbf{1}_{\{R^{-1}(t) \in F\}}\, dt = 0 \quad \text{a.s.}$$

By taking expectations, we see that

$$\int_0^\infty P(R^{-1}(t) \in F)\, dt = 0,$$

which proves the claim. □

The point is that the set $F'$ of $t$ such that $A$ fails at $R^{-1}(t)$ may be chosen deterministically. If $t \notin F'$, then, with probability one, $A_{R^{-1}(t)}$ holds, even though a priori we only knew this property for fixed, deterministic levels.

As a consequence of Lemma 26, we may choose a deterministic $t$ such that the limit theorem for the number of vertices on $\mathbf{T}^u$ remains true for the level $u = R^{-1}(t)$. For simplicity, we will assume that $t = 1$ is a valid choice, and we write $\mathcal{T}_0 = \mathbf{T}^{R^{-1}(1)}$ for the tree which has a set of vertices at level $t$ ($0 \le t \le 1$) given by the excursions above $R^{-1}(1)t$ that reach level $R^{-1}(1)$. Hence,

$$(26) \qquad \lim_{t \to 0} t^{1/(\alpha - 1)} |\mathcal{T}_0(1 - t)| \to K(R^{-1}(1))^{-1/\alpha - 1} Z_{R^{-1}(1)} \quad \text{a.s.}$$

The only thing that remains to be considered is the behavior of $t \mapsto R^{-1}(1 - t)$ when $t$ is small.

LEMMA 27. *As $t \to 0$, the following asymptotics hold almost surely:*

$$R^{-1}(1) - R^{-1}(1 - t) \sim t \frac{1}{\alpha(\alpha - 1)\Gamma(\alpha)} Z_{R^{-1}(1)}^{\alpha - 1},$$

*meaning that the ratio of the two sides converges to 1 almost surely.*

PROOF. Let

$$(27) \qquad\qquad\qquad q = \frac{1}{\alpha(\alpha - 1)\Gamma(\alpha)} Z_{R^{-1}(1)}^{\alpha - 1}.$$

The lemma follows simply from the fact that almost surely the function $R(t)$ is differentiable at $t = R^{-1}(1)$ since $Z$ is continuous at $R^{-1}(1)$. Its derivative is given by

$$\alpha(\alpha - 1)\Gamma(\alpha) Z_{R^{-1}(1)}^{1-\alpha} = q^{-1}$$

which is nonzero almost surely. Therefore, $R^{-1}(t)$ is also differentiable at $t=1$ and its derivative is $q$. □

PROOF OF THEOREM 4. Now, to finish the proof, note that for $t \leq 1$,
$$N(t) = \theta^{R^{-1}(1)}(R^{-1}(1-t)).$$
Since $R^{-1}(1-t) = R^{-1}(1) - tq + o(t)$, by monotonicity of $\theta^{R^{-1}(1)}$, we see that
$$N(t) \sim \theta^{R^{-1}(1)}(R^{-1}(1) - tq).$$
On the other hand, by (26), we have
$$\left(t\frac{q}{R^{-1}(1)}\right)^{1/(\alpha-1)} \theta^{R^{-1}(1)}(R^{-1}(1)-tq) \to K(R^{-1}(1))^{-1/(\alpha-1)} Z_{R^{-1}(1)} \qquad \text{a.s.}$$
After cancellation, we obtain that almost surely
$$t^{1/(\alpha-1)} N(t) \to K(\alpha(\alpha-1)\Gamma(\alpha))^{1/(\alpha-1)} = (\alpha\Gamma(\alpha))^{1/(\alpha-1)},$$
as stated in Theorem 4. □

5.3. *Evans' metric space and multifractal spectrum.* We begin with by a description of the basic setup for this section, which is Evans' random metric space $S$. This space was introduced by Evans in [23] in the case of Kingman's coalescent, and some properties of $S$ (such as its Hausdorff and packing dimensions) were derived in [8] in the case of a Beta$(2-\alpha, \alpha)$-coalescent and other coalescents behaving similarly (see [8], Theorem 1.7). The space $S$ is defined as the completion of $\mathbb{N}$ for the distance $d_S$ which is defined on $\mathbb{N}$ by
$$d_S(i,j) = \inf\{t : i \sim_{\Pi(t)} j\},$$
that is, $d_S(i,j)$ is the collision time of $i$ and $j$. Observe that $d_S$ is, in fact, an *ultrametric*, both on $\mathbb{N}$ and $S$, that is,
$$d_S(x,z) \leq d_S(x,y) \vee d_S(y,z) \qquad \forall x,y,z \in S.$$
The space $(S, d_S)$ is complete by definition and hence it is compact as soon as $\Pi(t)$ comes down from infinity. Indeed, for each $t > 0$, one needs only $N(t) < \infty$ balls of diameter $t$ to cover it, which implies that $S$ is precompact. Together with completeness, this makes the space $S$ compact. Given $B \subseteq S$, we write cl$B$ or $\bar{B}$ for its closure (with respect to $d_S$). Let $I_i(t) := \min\{j \in B_i(t)\}$ be the least element of $B_i(t)$. Then the set
$$\begin{aligned} U_i(t) &= \text{cl} B_i(t) \\ &= \text{cl}\{j \in \mathbb{N} : j \sim_{\Pi(t)} I_i(t)\} \\ &= \text{cl}\{j \in \mathbb{N} : d(j, I_i(t)) \leq t\} \\ &= \{y \in S : d(y, I_i(t)) \leq t\} \end{aligned}$$

is a closed ball with diameter at most $t$. The closed balls of $S$ are also the open balls of this space and every ball is of the form $U_i(t)$. In particular, it is easily seen that the collection of balls is countable. For $x \in S$ and $t \geq 0$, we write $\bar{B}_x(t)$ for the ball of center $x$ and diameter $t$ [observe that in the case $x \in \mathbb{N}$, this notation is consistent with the blocks convention for $\Pi(t)$].

It is possible (see [23]) to define almost surely a random measure $\eta(\cdot)$ on $S$ by requiring that for all $i \in \mathbb{N}$ and all $t \geq 0$, the measure $\eta(U_i(t))$ is the frequency of the block of $\Pi(t)$ containing $i$. We call $\eta$ the *mass-measure* or the *size-biased picking measure*. Recall that for $\gamma \leq 1/(\alpha - 1)$, the subset $S(\gamma)$ of $S$ is defined as

$$S(\gamma) = \left\{ x \in S : \liminf_{r \to 0} \frac{\log(\eta(\bar{B}_x(r)))}{\log r} \leq \gamma \right\}.$$

Results from [8] suggest that $\gamma = 1/(\alpha - 1)$ is the typical exponent for the size of a block as time goes down to 0. Hence, here, we are looking for existence of blocks whose size is abnormally large compared to the typical size as time goes down to 0. The next result gives the precise value of the Hausdorff dimension of this set (with respect to the distance on the space $S$).

The key idea for the proof of Theorem 5 is the observation that the space $S$, equipped with its mass measure $\eta$, can be thought of as the boundary of some Galton–Watson tree [more precisely, the reduced tree at level $R^{-1}(t)$] with the associated branching measure. Hence, the multifractal spectrum of $\eta$ in $S$ is the same as the multifractal spectrum of the branching measure in the boundary of a supercritical Galton–Watson tree. The case where the offspring distribution is heavy-tailed and has infinite variance has been recently studied by Mörters and Shieh [41] and we can use their result to conclude. For basic properties of the branching measure of a Galton–Watson tree, we recommend the references [37, 38, 41].

Recall that $\mathbf{T}^u$ designates the reduced tree at level $u$, that is, it is the tree where, for each level $0 \leq t \leq 1$, each vertex at level $t$ corresponds to one excursion of $H$ above level $ut$ that reaches level $u$. For our purposes, we eventually wish to work under the law of $(H(s), 0 \leq s \leq T_r)$ (conditionally on the event $\sup_{s \leq T_r} H_s > u$, otherwise the tree is empty), but it will sometimes be more convenient to use $N_u(\cdot)$, the excursion measure conditioned to hit level $u$. The difference is, of course that in the latter case, $\mathbf{T}^u$ is a tree with a single ancestor, while in the former case, $\mathbf{T}^u$ is actually a collection of a Poissonian number of i.i.d. trees joined at the root. These trees have the distribution of the reduced tree under $N_{(u)}(\cdot)$. We emphasize that for this study of the multifractal spectrum, this does not create any real difference.

By definition, a *ray* of $\mathbf{T}^u$ is a path $(\zeta(t), 0 \leq t \leq 1)$ such that $\zeta(0)$ is the root, for every $t, \zeta(t)$ is a vertex at level $t$ in $\mathbf{T}^u$ and for all $s \leq t$, $\zeta(s)$ is an ancestor of $\zeta(t)$. Then the boundary of the tree $\mathbf{T}^u$, denoted $\partial \mathbf{T}^u$, is just the set of all rays. The boundary $\partial \mathbf{T}^u$ can be equipped with a metric $\text{dist}_{\partial \mathbf{T}}$ by letting $\text{dist}_{\partial \mathbf{T}}(U, V) = 1 - t$ if $t$ is the height at which $U$ and $V$ diverge. Let $|\mathbf{T}^u(t)| := \theta^{(u)}(ut)$ be the size of generation at level $t$. By Proposition 19, we see that $(|\mathbf{T}^u(1 - e^{-t})|, t \geq 0)$ is a continuous-time Galton–Watson process where individuals live for an exponential time with parameter 1 and then reproduce with offspring distribution $\chi$. Recall from Lemma 21 that there is a random variable $W > 0$ almost surely such that

$$W = \lim_{t \to \infty} e^{-t/(\alpha-1)} |\mathbf{T}^u(1 - e^{-t})|.$$

Furthermore, for every vertex $v \in \mathbf{T}^u$, we can define $\mathbf{T}^u(v)$, the subtree rooted at $v$, and $W(v)$, the limit (which exists almost surely) of its associated martingale. As there are countably many branching points of $\mathbf{T}^u$, this allows one to build a natural measure $\mu$, called the *branching measure* on $\partial \mathbf{T}^u$, by introducing the requirement

(28) $$\mu(\{\zeta \in \partial \mathbf{T}^u : \zeta(1 - e^{-t}) = v\}) = \frac{W(v)}{e^{t/(\alpha-1)}}.$$

Observe that the set on the left-hand side is a ball of radius $e^{-t}$ centered on any ray $\zeta$ such that $\zeta(1 - e^{-t}) = v$. Having defined $\mu$ on arbitrary balls of the boundary of the tree $\partial \mathbf{T}^u$, this uniquely extends to a measure $\mu$ which is defined on arbitrary subsets of $\partial \mathbf{T}^u$ by Carathéodory's Extension theorem (see page 438 of [19]).

When $u > 0$ is a fixed deterministic level, $\mathbf{T}^u$ is a collection of Galton–Watson trees. The definitions introduced above then coincide with the standard notions of distance, boundary and branching measure for a collection of Galton–Watson trees. The lemma below is essentially a reformulation of Theorems 2.1 and 2.2 in [41] within our framework.

LEMMA 28. *Conditionally on* $\sup_{0 \leq s \leq T_r} H_s > u$, *the multifractal spectrum of $\mu$ is given as follows: for all* $\frac{1}{\alpha} \leq \gamma \leq \frac{1}{\alpha - 1}$,

$$\dim_{\mathcal{H}} \left\{ V \in \partial \mathbf{T}^u : \liminf_{r \to 0} \frac{\log(\mu(B(V, r)))}{\log r} \leq \gamma \right\} = \gamma \alpha - 1$$

*and the set is empty if* $\gamma < 1/\alpha$. *If* $\frac{1}{\alpha - 1} < \gamma \leq \frac{\alpha}{(\alpha-1)^2}$, *then*

$$\dim_{\mathcal{H}} \left\{ V \in \partial \mathbf{T} : \limsup_{r \to 0} \frac{\log(\mu(B(V, r)))}{\log r} \geq \gamma \right\} = \frac{\alpha}{\gamma(\alpha - 1)^2} - 1$$

*and the set is empty when* $\gamma > \alpha/(\alpha - 1)^2$.

PROOF. First, we remark that it suffices to prove this result under the measure $N_{(u)}$. Moreover, it is elementary to check that

$$\left\{V \in \partial \mathbf{T}^u : \liminf_{t \to \infty} \frac{\log(\mu(B(V, e^{-t})))}{-t} = \gamma\right\}$$
$$= \left\{V \in \partial \mathbf{T}^u : \liminf_{n \to \infty} \frac{\log(\mu(B(V, e^{-n})))}{-n} = \gamma\right\}$$

and that

$$\left\{V \in \partial \mathbf{T}^u : \limsup_{t \to \infty} \frac{\log(\mu(B(V, e^{-t})))}{-t} = \gamma\right\}$$
$$= \left\{V \in \partial \mathbf{T}^u : \limsup_{n \to \infty} \frac{\log(\mu(B(V, e^{-n})))}{-n} = \gamma\right\}.$$

Sampling at these discrete times gives us a discrete-time Galton–Watson process which satisfies the assumptions of Theorems 2.1 and 2.2 of [41]. Its offspring variable is given by

$$\chi_{\text{discrete}} := |\mathbf{T}^u(1 - e^{-1})|.$$

Observe that, by construction, $P(\chi_{\text{discrete}} = 0) = 0$ and $P(\chi_{\text{discrete}} = 1) < 1$. Furthermore, it is easily seen that $E(\chi_{\text{discrete}}) = e^{1/(\alpha-1)}$. By [5], Corollary 2, Chapter III.6, the offspring variable $\chi_{\text{discrete}}$ in discrete time and $\chi$ satisfy the $X \log X$ condition simultaneously, so

$$E(\chi_{\text{discrete}} \log \chi_{\text{discrete}}) < \infty.$$

The last step is to check the values of the two constants

$$\tau := -\log(P(\chi_{\text{discrete}} = 1))/\log(E(\chi_{\text{discrete}}))$$

and

$$r := \liminf_{x \to \infty} \frac{-\log P(\chi_{\text{discrete}} > x)}{\log x}.$$

Note that $\chi_{\text{discrete}} = 1$ occurs if the ancestor has not reproduced by time 1. Since the time at which she reproduces is, on this timescale, an exponential random variable with mean 1, we see that $P(\chi_{\text{discrete}} = 1) = e^{-1}$, so $\tau = -(\alpha - 1)\log(e^{-1}) = \alpha - 1$. To compute $r$ requires a few more arguments. Now, it is known (see [37], (3.1) and (3.2)) that $r$ is equal to

$$\sup\{a > 0, E(\chi_{\text{discrete}}^a) < \infty\}.$$

On the other hand, by [5], Corollary 1, Chapter III.6 for all $a > 1$, $E(\chi_{\text{discrete}}^a) < \infty$ if and only if $E(\chi^a) < \infty$. Using (10), we see that $\chi$ admits moments of order up to and excluding $\alpha$, therefore $r = \alpha$. Application of Theorems 2.1 and 2.2 of [41] concludes the proof of the lemma. □

The proof of Theorem 5 is now straightforward. We show that the multifractal spectrum of $\eta$ in $S$ with respect to the metric $d_S$ is necessarily the same as the multifractal spectrum of $\mu$ in $\partial \mathbf{T}$ with respect to $\text{dist}_{\partial \mathbf{T}}$.

PROOF OF THEOREM 5. Let $\mathcal{T}$ be the tree whose vertices at level $t$ consist of those excursions above level $R^{-1}(t)$ that reach level $R^{-1}(1)$. As above, the boundary $\partial \mathcal{T}$ of the tree $\mathcal{T}$ is just the set of all "infinite" paths, that is, of paths $(\zeta(t), 0 \leq t \leq 1)$ such that for every $t$, $\zeta(t)$ is at level $t$ of $\mathcal{T}$. We may equip $\partial \mathcal{T}$ with the following metric: the distance between two rays $\zeta$ and $\zeta'$ is simply

$$\text{dist}_{\partial \mathcal{T}}(\zeta, \zeta') = 1 - \sup\{t \leq 1 : \zeta(t) = \zeta'(t)\}.$$

There is a one-to-one map $\Phi$ between $S$ and $\partial \mathcal{T}$ which can be described as follows: let $\zeta \in \partial \mathcal{T}$, then for each $t \in (0,1)$, the vertex $\zeta(1-t)$ corresponds, by definition, to an excursion above $R^{-1}(1-t)$ that hits level $R^{-1}(1)$ and hence to a block $B_\zeta(t)$ of the partition $\Pi(t)$, where $\Pi$ is the embedded coalescent process. When $t < t'$, $B_\zeta(t) \subseteq B_\zeta(t')$. Define $i(t) := \min B_\zeta(t)$, the least element of the block that corresponds to the vertex $\zeta(1-t)$. Note that the function $i(t)$ satisfies the Cauchy criterion (with respect to the metric $d_S$) as $t \to 0$, by construction. Since $S$ is a complete metric space under $d_S$, it follows that there is a unique $x \in S$ such that $d_S(x, i(t)) \to 0$ when $t \to 0$. We put $\Phi^{-1}(\zeta) = x$. In the converse direction, since $\mathbb{N}$ is dense in $S$, for any $x \in S$, we may consider a sequence $(i_n, n = 1, 2, \ldots)$ in $\mathbb{N}$ such that $d_S(i_n, x) \to 0$ when $n \to \infty$. Without loss of generality, we may assume that $d_S(i_n, x)$ is monotone decreasing. Then the sequence of blocks $B(i_n, t_n)$ that contain $i_n$ at time $t_n = d_S(i_n, x)$ defines a unique ray $\zeta_x$ such that $\zeta_x(1 - t_n)$ corresponds to $B(i_n, t_n)$ for each $n$. Moreover, $\zeta_x$ does not depend on the particular sequence $i_n$ converging to $x$, so we may unambiguously define $\Phi(x) = \zeta_x$. It is easy to see that $\Phi(\Phi^{-1}(\zeta)) = \zeta$. For instance, this map $\Phi$ acts on the integers as follows: $\forall i \in \mathbb{N}$, $\Phi(i)$ is the unique ray $(\zeta(t), t \geq 0)$ such that for each $t \geq 0$, the integer $i$ is in the block of $\Pi(t)$ which corresponds to $\zeta(t)$.

Hence, we may identify $S$ with $\partial \mathcal{T}$ and note that, by construction, distances are preserved in this identification

$$d_S(x, y) = \text{dist}_{\partial \mathcal{T}}(\Phi(x), \Phi(y)).$$

Furthermore, if $z$ is a vertex at level $t$ of $\mathcal{T}$, let $\ell(z)$ the total local time at level $R^{-1}(1)$ of the excursion defining $z$, divided by $Z_{R^{-1}(1)}$, the total local time of the whole process $(H_s, s \leq T_r)$ at level $R^{-1}(1)$. The correspondence between local time at level $R^{-1}(1)$ and asymptotic frequencies of the blocks of $\Pi(t)$ implies that

$$\eta(B(x, t)) = \ell(\zeta_x(t)), \tag{29}$$

where $\zeta_x(t)$ is the vertex corresponding to $B(x,t)$, that is, the vertex at level $t$ on the ray $\zeta_x$. Hence, as the map $\Phi$ preserves the distance, it is easy to see that $\dim_{\mathcal{H}} S(\gamma) = \dim_{\mathcal{H}} S'(\gamma)$, where

$$S'(\gamma) = \left\{ \zeta \in \partial \mathcal{T} : \liminf_{t \to 0} \frac{\log \ell(\zeta(1-t))}{\log t} \leq \gamma \right\} \tag{30}$$

because the two sets coincide via the map $\Phi$. Thus, we want to prove that $\dim_{\mathcal{H}} S'(\gamma) = (\gamma \alpha - 1)$. On the other hand, recall that $\mathcal{T}$ is just a rescaling of $\mathcal{T}_0$, which is the shorthand notation for the reduced tree at level $R^{-1}(1)$. Recall that this tree has a set of vertices at level $t$ (for $0 \leq t \leq 1$) corresponding to excursions above level $tR^{-1}(1)$ reaching $R^{-1}(1)$. Let us first treat the case $\gamma \leq 1/(\alpha - 1)$ of "thick points." By Lemma 27, we have that $R^{-1}(1) - R^{-1}(1-t) \sim tq$ (where, as before, $q$ denotes the random number in Lemma 27), so it is enough to prove that $\dim_{\mathcal{H}} S'_0(\gamma) = (\gamma \alpha - 1)$, where

$$S'_0(\gamma) = \left\{ \zeta \in \partial \mathcal{T}_0 : \liminf_{t \to 0} \frac{\log \ell(\zeta(1-t))}{\log t} \leq \gamma \right\}.$$

On the other hand, by Lemma 23, for a fixed level $u > 0$, the limit $W$ of the Kesten–Stigum martingale associated with the reduced tree at level $u$ is a constant multiple of the local time at $u$. Let $(\zeta(t), 0 \leq t \leq 1)$ be a ray in $\partial \mathbf{T}^u$. Applying this to the subtree rooted at $v = \zeta(t)$, it follows that the number $W(v)$ defining the branching measure on $\partial \mathbf{T}^u$ is also a constant multiple of the local time $\ell(v)$ at level $u$ enclosed in the excursion corresponding to vertex $v$,

$$W(v) = K(ue^{-t})^{-1/(\alpha-1)} \ell(\zeta(t)),$$

since the subtree rooted at $v$ has the law of $\mathbf{T}^{ue^{-t}}$. In other words, dividing both sides by $e^{t/(\alpha-1)}$ and referring to (28), if $\mu$ is the branching measure on $\partial \mathbf{T}^u$, then almost surely, for all $t > 0$,

$$\mu(B(\zeta, e^{-t})) = K\ell(\zeta(t)),$$

that is, the branching measure associated with a vertex $\zeta(t) = z \in \mathbf{T}^u$ is a constant multiple of the local time $\ell(z)$ enclosed at level $u$ in the excursion corresponding to $z$. Therefore, using Lemma 28, this implies that almost surely, conditionally on the event $\sup_{0 \leq s \leq T_r} H_s > u$,

$$\dim_{\mathcal{H}} \left\{ \zeta \in \mathbf{T}^u : \liminf_{t \to 0} \frac{\log \ell(\zeta(1-t))}{\log t} \leq \gamma \right\} = \gamma \alpha - 1.$$

We may therefore apply Lemma 26 to conclude that if $t \notin N$, where $N$ is a deterministic set of Lebesgue measure zero, then this property also holds for the reduced tree at level $R^{-1}(t)$. There is, of course, no loss of generality in assuming that $1 \notin N$, so we conclude that

$$\dim_{\mathcal{H}} S'_0(\gamma) = \gamma \alpha - 1,$$

as required. When $\gamma > 1/(\alpha - 1)$, the proof follows the same lines and uses the "thin points" part of Lemma 28. This concludes the proof of Theorem 5. □

**6. Site and allele frequency spectrum.** Our goal in this section is to prove Theorem 9. Our proof relies heavily on the connection between Beta-coalescents and Galton–Watson processes developed in the previous section. Throughout this section, $(\xi_t, t \geq 0)$ will denote the continuous-time Galton–Watson process where individuals live for an independent exponential amount of time and then give birth to a number of offspring distributed according to $\chi$, where $P(\chi = 0) = P(\chi = 1) = 0$ and, for $k \geq 2$,

$$P(\chi = k) = \frac{\alpha(2-\alpha)(3-\alpha)\cdots(k-1-\alpha)}{k!} = \frac{\alpha \Gamma(k-\alpha)}{k! \Gamma(2-\alpha)}.$$

This offspring distribution is supercritical, with mean $m = 1 + 1/(\alpha - 1)$. Also, recall that $M_k(n)$ denotes the number of families of size $k$ in the infinite sites model when the sample has $n$ individuals, and that $N_k(n)$ is the equivalent quantity in the infinite alleles model.

6.1. *Expected values.* Suppose marks occur at times of a constant rate $\theta$ Poisson point process along the branches of a reduced tree at level 1 under the measure $N_{(1)}$, so that the reduced tree has a single ancestor. Recall that the number of branches of $\mathbf{T}^1$ at level $1 - e^{-t}$ is a Galton–Watson process. Hence, after rescaling, this amounts to having mutation marks at intensity $\theta e^{-s}$ per unit length at time $s$ on the Galton–Watson tree that comes from the process $\xi$. We will stop the Galton–Watson process at a fixed time $t$. If there is a mutation at time $s < t$, then we say that it creates a family of size $k$ if the individual with the mutation at time $s$ has $k$ descendants alive in the population at time $t$. Let $M_k^{\mathrm{GW}}(t)$ denote the number of families of size $k$ at time $t$. The following result shows that a simple calculation gives the asymptotic behavior of $E[M_k^{\mathrm{GW}}(t)]$. A sharper argument will be needed to establish convergence in probability.

PROPOSITION 29. *Let $\tau$ be an independent exponential random variable with mean $1/c$, where*

$$c = \frac{2-\alpha}{\alpha - 1}.$$

*We have*

$$\lim_{t \to \infty} e^{-ct} E(M_k^{\mathrm{GW}}(t)) = \frac{\theta}{c} P(\xi_\tau = k).$$

PROOF. By applying the branching property and using the facts that $E[\xi_t] = e^{(m-1)t}$ and that $m - 2 = c$ for the third equality, we obtain

$$\begin{aligned} E(M_k^{\mathrm{GW}}(t)) &= \int_0^t P(\text{there is a mark in } dl) P(\xi_{t-l} = k) \\ &= \int_0^t E(\xi_l \theta e^{-l}) \, dl \ P(\xi_{t-l} = k) \\ &= \int_0^t \theta e^{cl} P(\xi_{t-l} = k) \, dl \\ &= \theta e^{ct} \int_0^t e^{-cu} P(\xi_u = k) \, du \\ &= e^{ct} \frac{\theta}{c} \int_0^t c e^{-cu} P(\xi_u = k) \, du. \end{aligned}$$

Multiplying both sides by $e^{-ct}$ and letting $t \to \infty$, we get

$$\begin{aligned} \lim_{t \to \infty} e^{-ct} E(M_k^{\mathrm{GW}}(t)) &= \frac{\theta}{c} \int_0^\infty c e^{-cu} P(\xi_u = k) \, du \\ &= \frac{\theta}{c} P(\xi_\tau = k). \end{aligned} \qquad \square$$

To make the limiting expression for $E[M_k^{\mathrm{GW}}(t)]$ more explicit, we now calculate $P(\xi_\tau = k)$.

LEMMA 30. *For all positive integers $k$, we have*

$$(31) \qquad P(\xi_\tau = k) = \frac{(2-\alpha)\Gamma(k+\alpha-2)}{\Gamma(\alpha-1)k!}.$$

PROOF. We prove the result by induction. Note that $\xi_\tau = 1$ if and only if there are no birth events before time $\tau$. Because $\tau$ has an exponential distribution with rate parameter $c$ and individuals give birth at rate 1, it follows that

$$P(\xi_\tau = 1) = \frac{c}{1+c} = 2 - \alpha,$$

which agrees with the right-hand side of (31) when $k = 1$.

Now, suppose that $k \geq 2$ and (31) is valid for $j = 1, \ldots, k-1$. Let $r_k = P(\xi_t = k$ for some $t \leq \tau)$. By conditioning on the number of individuals before there were $k$ individuals, we get

$$(32) \qquad r_k = \sum_{j=1}^{k-1} r_j \cdot \frac{j}{j+c} P(\chi = k - j + 1)$$

because if there are $j < k$ individuals, then the probability of having another birth before time $\tau$ is $j/(j+c)$ and if this happens, the probability that there are $k$ individuals after the next birth is $P(\chi = k - j + 1)$. If $\xi_t = k$ for some $t \leq \tau$, then we will have $\xi_\tau = k$ if and only if $\tau$ occurs before the next birth event. When there are $k$ individuals, birth events happen at rate $k$, so the probability that $\tau$ happens before the next birth is $c/(k+c)$. Therefore, $P(\xi_\tau = k) = cr_k/(k+c)$ and so $r_k = P(\xi_\tau = k)(k+c)/c$. Substituting this into (32), we get

$$(33) \qquad P(\xi_\tau = k) = \frac{1}{k+c} \sum_{j=1}^{k-1} j P(\xi_\tau = j) P(\chi = k - j + 1).$$

Using the induction hypothesis and the fact that $P(\chi = k) = \alpha \Gamma(k - \alpha)/(k! \Gamma(2 - \alpha))$, we obtain

$$P(\xi_\tau = k) = \frac{\alpha(2-\alpha)}{\Gamma(\alpha - 1)\Gamma(2-\alpha)(k+c)} \sum_{j=1}^{k-1} \frac{\Gamma(j+\alpha-2)\Gamma(k-j+1-\alpha)}{(j-1)!(k-j+1)!}.$$

Using the fact that $k + c = (k\alpha - k + 2 - \alpha)/(\alpha - 1)$ and letting $\ell = j - 1$ in the sum, we get

$$(34) \qquad \begin{aligned} P(\xi_\tau = k) &= \frac{\alpha(\alpha-1)(2-\alpha)}{(k\alpha - k + 2 - \alpha)\Gamma(\alpha-1)\Gamma(2-\alpha)} \\ &\quad \times \sum_{\ell=0}^{k-2} \frac{\Gamma(\ell+\alpha-1)\Gamma(k-\ell-\alpha)}{\ell!(k-\ell)!}. \end{aligned}$$

If $a, b \in \mathbb{R}$ and $n \in \mathbb{N}$, then by starting with the identity $(1-x)^{-a}(1-x)^{-b} = (1-x)^{-(a+b)}$ and considering the $n$th order term in the Taylor series expansion of both sides, we get (see, e.g., page 70 in [3])

$$\sum_{k=0}^{n} \frac{(a)_k (b)_{n-k}}{k!(n-k)!} = \frac{(a+b)_n}{n!},$$

where $(a)_k = a(a+1)\cdots(a+k-1)$. Since $(a)_k = \Gamma(a+k)/\Gamma(a)$, it follows that

$$(35) \qquad \sum_{k=0}^{n} \frac{\Gamma(a+k)\Gamma(b+n-k)}{k!(n-k)!} = \frac{\Gamma(a)\Gamma(b)(a+b)_n}{n!}.$$

When $a + b = -1$, we have $(a+b)_n = 0$. Therefore, (35) with $a = \alpha - 1$ and $b = -\alpha$ implies that the sum on the right-hand side of (34) would be zero if it went up to $k$ rather than $k - 2$. It follows that the sum up to $k - 2$ is equal to the negative of the sum of the terms when $\ell = k$ and $\ell = k-1$, which is

$$-\frac{\Gamma(k+\alpha-2)\Gamma(1-\alpha)}{(k-1)!} - \frac{\Gamma(k-\alpha-1)\Gamma(-\alpha)}{k!}$$

$$= \frac{\Gamma(2-\alpha)\Gamma(k+\alpha-2)(k\alpha-k+2-\alpha)}{k!\alpha(\alpha-1)}.$$

Combining this result with (34) gives (31). The lemma follows by induction. □

6.2. *A queueing system result.* The problem on a Galton–Watson tree will essentially reduce to the following lemma.

Let $Q_t$ be the state of a queueing system where customers arrive at rate $Ae^{ct}$ for some constants $A$ and $c > 0$. We assume that there are infinitely many servers and that each customer requires an independent exponential rate $\lambda$ amount of time to be served, so when the state of the queue is $m$, the departure rate is $\lambda m$ per unit of time.

LEMMA 31. *As $t \to \infty$, almost surely*

$$e^{-ct}Q_t \to \frac{A}{\lambda+c}.$$

PROOF. Because all customers depart at rate $\lambda$, the number of customers at time zero does not affect the limiting behavior of the queue as $t \to \infty$. Therefore, we may assume that the number of customers at time zero is Poisson with mean $A/(\lambda+c)$. The probability that a customer who arrives at time $s \leq t$ is still in the queue at time $t$ is $e^{-\lambda(t-s)}$. Therefore, the distribution of $Q_t$ is Poisson with mean

$$\frac{Ae^{-\lambda t}}{\lambda+c} + \int_0^t Ae^{cs}e^{-\lambda(t-s)}\,ds = \frac{Ae^{ct}}{\lambda+c}.$$

For all positive integers $n$, let $t_n = (3/c)\log n$, so $E[Q_{t_n}] = An^3/(\lambda+c)$. Let $B_n$ be the event that $(1-\varepsilon)An^3/(\lambda+c) \leq Q_{t_n} \leq (1+\varepsilon)An^3/(\lambda+c)$. Note that if $Z$ has a Poisson distribution with mean $\mu$, then

(36) $$P(|Z-\mu| > \varepsilon\mu) \leq \frac{1}{\varepsilon^2\mu},$$

by Chebyshev's inequality. Applying (36) with $\mu = An^3/(\lambda+c)$, we get

$$P(B_n^c) \leq \frac{\lambda+c}{A\varepsilon^2 n^3}.$$

Therefore, by the Borel–Cantelli lemma, almost surely $B_n$ occurs for all but finitely many $n$.

Between times $t_n$ and $t_{n+1}$, the number of arrivals is Poisson with mean at most
$$\int_{t_n}^{t_{n+1}} Ae^{cs}\, ds = \frac{A}{c}((n+1)^3 - n^3) \leq \frac{3A(n+1)^2}{c}.$$
Therefore, the probability that there are more than $6A(n+1)^2/c$ arrivals between times $t_{n-1}$ and $t_n$ is at most the probability that a Poisson random variable with mean $3A(n+1)^2/c$ is greater than $6A(n+1)^2/c$, which by (36) with $\varepsilon = 1$, is at most $c/(3A(n+1)^2)$. The number of departures between times $t_n$ and $t_{n+1}$ also has a Poisson distribution and since $E[Q_t]$ is an increasing function of $t$, the expected number of departures between times $t_n$ and $t_{n+1}$ is also bounded by $3A(n+1)^2/c$. Therefore, the probability that there are more than $6A(n+1)^2/c$ departures between times $t_n$ and $t_{n+1}$ is at most $c/(3A(n+1)^2)$. Let $D_n$ be the event that between times $t_{n-1}$ and $t_n$, there are at most $6A(n+1)^2/c$ arrivals and at most $6A(n+1)^2/c$ departures. By the Borel–Cantelli lemma, almost surely $D_n$ occurs for all but finitely many $n$.

Suppose that $B_n$ and $D_n$ occur for all $n \geq N$. Suppose that $t_n \leq t \leq t_{n+1}$. If $n \geq N$, then
$$\frac{(1-\varepsilon)An^3}{\lambda+c} - \frac{6A(n+1)^2}{c} \leq Q_t \leq \frac{(1+\varepsilon)An^3}{\lambda+c} + \frac{6A(n+1)^2}{c}.$$
Because $1/n^3 \leq e^{-ct} \leq 1/(n+1)^3$, it follows that
$$\limsup_{n\to\infty} e^{-ct}Q_t \leq \frac{(1+\varepsilon)A}{\lambda+c} \qquad \text{a.s.}$$
and
$$\liminf_{n\to\infty} e^{-ct}Q_t \geq \frac{(1-\varepsilon)A}{\lambda+c} \qquad \text{a.s.}$$
Since $\varepsilon > 0$ is arbitrary, the result follows. □

Having proven this result, we easily deduce the following one. Suppose that $(Q_t, t \geq 0)$ is the length of the queue in a queueing system, where the arrival rate is a random process $a_t$ [i.e., the process of arrivals $(Q_t^+, t \geq 0)$ is a counting process such that $Q_t^+ - \int_0^t a_s\, ds$ is a martingale] and the departure rate at time $t$, which is nonrandom, is $\lambda(t)$ per customer. Then if $a_t$ and $\lambda(t)$ have the correct asymptotics as $t \to \infty$, the asymptotics of $Q_t$ are also the same as in the previous case.

LEMMA 32. *If $a_t \sim Ae^{ct}$ almost surely as $t \to \infty$ and $\lim_{t\to\infty} \lambda(t) = \lambda$, then almost surely*
$$e^{-ct}Q_t \to \frac{A}{\lambda+c}.$$

PROOF. Let $A_t = \int_0^t a_s\, ds$. Since $Q_t^+$ is a counting process and $Q_t^+ - A_t$ is a martingale, there exists a Poisson process $N_t^0$ such that $Q_t^+ = N_{A_t}^0$. Let $\varepsilon > 0$ and consider the function

$$b_t = (A(1+\varepsilon)e^{ct} - a_t)_+.$$

Let $N^1$ be an independent Poisson process. Compare the state of the queue $(Q_t, t \geq 0)$ with the queue $(Q_{1,t}, t \geq 0)$, in which customers arrive with the jumps of $(N_{A_t}^0 + N_{B_t}^1, t \geq 0)$, where $B_t = \int_0^t b_s\, ds$, and customers get served at rate $\lambda(t)$. By properties of Poisson processes, the arrival process of the queue $Q_1$ is thus itself a Poisson process with rate $a_t + b_t$ per unit time. Observe that for $t$ sufficiently large, $a_t \leq A(1+\varepsilon)e^{ct}$, so for $t$ sufficiently large, $b_t = A(1+\varepsilon)e^{ct} - a_t$. Thus, for $t$ sufficiently large, the total rate of arrivals for the queue $Q_1$ is $a_t + b_t = A(1+\varepsilon)e^{ct}$. Let $(Q_{2,t}, t \geq 0)$ be the queue where arrivals are given by $N_{A_t}^0 + N_{B_t}^1$ when $a_t \leq A(1+\varepsilon)e^{ct}$ and $N_{A(1+\varepsilon)e^{ct}}^2$ otherwise, where $(N_t^2, t \geq 0)$ is another independent Poisson process. Again, assume that customers depart from the queue at rate $\lambda(t)$. Since customers depart from $Q_1$ and $Q_2$ at the same rate, the queues can be coupled so that they are identical after a certain random time $T$. Moreover, $(Q_{2,t}, t \geq 0)$ is a queueing system where arrivals occur at rate $A(1+\varepsilon)e^{ct}$ throughout time. Because $\lambda(t) \to \lambda$, we have $\lambda(t) \geq \lambda - \varepsilon$ for sufficiently large $t$. Therefore, the queue $(Q_{2,t}, t \geq 0)$ can be coupled with another queue $(Q_{3,t}, t \geq 0)$ with arrival rate $A(1+\varepsilon)e^{ct}$ and departure rate $\lambda - \varepsilon$ such that $Q_{2,t} \leq Q_{3,t}$ for sufficiently large $t$, because for $t$ sufficiently large, all customers depart $Q_2$ at least as quickly as they depart $Q_3$. Hence, by Lemma 31, almost surely

$$\limsup_{t \to \infty} e^{-ct} Q_{2,t} \leq \frac{A(1+\varepsilon)}{(\lambda - \varepsilon) + c}$$

and similarly for $Q_1$, because $Q_1$ and $Q_2$ have the same asymptotics. By construction, we also have that for all $t \geq 0$, $Q_t \leq Q_{1,t}$, because every customer who arrives in $Q$ also arrives in $Q_1$. By taking $\varepsilon \to 0$, this implies that

$$\limsup_{t \to \infty} e^{-ct} Q_t \leq \frac{A}{\lambda + c}.$$

Applying similar reasoning, we get $\liminf_{t \to \infty} e^{-ct} Q_t \geq (1-\varepsilon)A/((\lambda+\varepsilon)+c)$ a.s. for all $\varepsilon > 0$, which implies the lemma. □

6.3. *Almost sure result for a Galton–Watson tree.* Recall that we are considering the Galton–Watson tree associated with the branching process $(\xi_t, t \geq 0)$, with mutation marks along the branches at rate $\theta e^{-s}$ at time $s$. By Lemma 21, there is a random variable $W$ such that

$$e^{-(m-1)t} \xi_t \to W \qquad \text{a.s.}$$

Recall that $M_k^{\mathrm{GW}}(t)$ denotes the number of marks before time $t$ such that the individual who gets the mutation has $k$ descendants at time $t$. Likewise, let $N_k^{\mathrm{GW}}(t)$ denote the number of blocks of size $k$ in the allelic partition at time $t$, when we assume two individuals have different alleles if any of their ancestors have had a mutation since their most recent common ancestor.

For the proof, we introduce two other quantities. Let $L_k(t)$ be the number of mutations before time $t$ such that the individual who gets the mutation has $k$ descendants alive at time $t$ and none of this individual's descendants undergoes another mutation before time $t$. Let $K(t)$ be the number of mutations before time $t$ such that some descendant of the individual that undergoes the mutation also undergoes another mutation before time $t$. The strategy of the proof will be to show that $M_k^{\mathrm{GW}}(t)$ and $N_k^{\mathrm{GW}}(t)$ both behave asymptotically like $L_k(t)$, while $K(t)$ is of lower order. The lemma below concerns $L_k(t)$.

LEMMA 33. *For all $k \geq 1$,*
$$e^{-ct} L_k(t) \to \frac{\theta W}{c} P(\xi_\tau = k) \qquad a.s.$$

PROOF. Our first step is to prove that this result holds with a limit being $\theta W a_k$ for some deterministic sequence of positive numbers $a_k$.

We prove this by induction on $k \geq 1$. For $k = 1$, observe that, conditionally on the process $(\xi_t, t \geq 0)$, the process $L_1(t)$ can be viewed as a birth-and-death chain in which the total birth rate is $\theta e^{-t} \xi_t$ and each individual dies at rate $1 + \theta e^{-t}$. Indeed, $L_1(t)$ increases by one every time some branch gets hit by a mutation. Since marks arrive at rate $\theta e^{-t} \, dt$ at time $t$ on each branch of the Galton–Watson tree, this means that, conditional on $(\xi_t, t \geq 0)$, new mutations occur at rate $\theta e^{-t} \xi_t \, dt$. Also, $L_1(t)$ decreases by one each time a member of a family of size 1 either reproduces or experiences a mutation, which happens at rate $1 + \theta e^{-t}$ for every individual. Because $e^{-(m-1)t} \xi_t \to W$ a.s., we can view $L_1(t)$ as a queueing system whose arrival rate is asymptotic to $\theta W e^{ct}$ and whose departure rate converges to 1. Therefore, by conditioning on $W$ and applying Lemma 32, we have
$$e^{-ct} L_1(t) \to \frac{\theta W}{1+c}.$$

Because
$$P(\xi_\tau = 1) = \int_0^\infty c e^{-cu} P(\xi_u = 1) \, du$$
$$= \int_0^\infty c e^{-cu} e^{-u} \, du$$
$$= c/(c+1),$$

we can take $a_1 = (1/c)P(\xi_\tau = 1)$, which is deterministic.

Now suppose that $k \geq 2$. Note that families of size $k$ are obtained when an individual in a family of size $j$ with $j \leq k-1$ reproduces and has $k-j+1$ offspring. Therefore, the process $(L_k(t), t \geq 0)$ is a birth-and-death chain with arrival rate

$$\text{(37)} \qquad \sum_{j=1}^{k-1} jL_j(t)P(\chi = k-j+1)\, dt.$$

We emphasize that this does not mean that conditionally on $(L_j(t), t \geq 0, i = 1, \ldots, k-1)$, the process $L_k$ is a queueing system with arrival rate (37). Indeed, the positive jump times of $L_k$ are necessarily negative jump times of $L_j$ for some $j < k$. Instead, this means that the arrival process $L_k^+$ for the queue $L_k$ is a counting process such that

$$\text{(38)} \qquad L_k^+(t) - \int_0^t \sum_{j=1}^{k-1} jL_j(s)P(\chi = k-j+1)\, ds$$

is a martingale and conditionally on $L_k^+$, the process $L_k(t)$ is independent of the lower-level queues $L_j$, $j = 1, \ldots, k-1$. The departure rate at time $t$ is $k(1 + \theta e^{-t})$ because for each family of size $k$, there are $k$ individuals that could reproduce or experience mutation.

In particular, the arrival rate (37) for $L_k(t)$ is almost surely asymptotic to

$$\theta W \left( \sum_{j=1}^{k-1} ja_j P(\chi = k-j+1) \right) e^{ct}.$$

Applying Lemma 32 with $\lambda = k$, we conclude

$$e^{-ct} L_k(t) \to \theta W a_k \qquad \text{a.s.,}$$

where

$$a_k = \frac{1}{k+c} \sum_{j=1}^{k-1} ja_j P(\chi = k-j+1).$$

Thus, the constants $a_k$ satisfy the same recursion established in (33) for $P(\xi_\tau = k)$. Because $a_1 = (1/c)P(\xi_\tau = k)$, it follows that $a_k = (1/c)P(\xi_\tau = k)$ for all $k$. $\square$

We now use this result to obtain the asymptotic behavior of the quantities $M_k^{\text{GW}}$ and $N_k^{\text{GW}}$.

LEMMA 34. *For all $k \geq 1$, almost surely*

$$e^{-ct} M_k^{\text{GW}}(t) \to \frac{\theta W}{c} P(\xi_\tau = k)$$

*and*

$$e^{-ct} N_k^{\text{GW}}(t) \to \frac{\theta W}{c} P(\xi_\tau = k).$$

PROOF. Note that every mutation before time $t$ that is counted by $L_k(t)$ is inherited by $k$ individuals at time $t$. By the definition of $L_k(t)$, these $k$ individuals experience no additional mutations, so they form a block of the allelic partition at time $t$. It follows that $L_k(t) \leq M_k^{\text{GW}}(t)$ and $L_k(t) \leq N_k^{\text{GW}}(t)$. Furthermore, if any mutation not counted by $L_k(t)$ is passed on to $k$ individuals at time $t$ or gives rise to a block of size $k$ in the allelic partition at time $t$, then some descendant of the individual that experiences the mutation must experience another mutation before time $t$. Therefore, we have $M_k^{\text{GW}}(t) \leq L_k(t) + K(t)$ and $N_k^{\text{GW}}(t) \leq L_k(t) + K(t)$. Thus, the result will follow from Lemma 33 once we prove that

$$\lim_{t \to \infty} e^{-ct} K(t) = 0 \quad \text{a.s.} \tag{39}$$

To prove (39), note that if $M(t)$ denotes the total number of mutations before time $t$, then for all positive integers $N$,

$$K(t) = M(t) - \sum_{k=1}^{\infty} L_k(t)$$

$$\leq M(t) - \sum_{k=1}^{N} L_k(t).$$

Conditional on $(\xi_t, t \geq 0)$, the process $(M(t), t \geq 0)$ is a queueing system with departure rate zero and arrival rate $\theta e^{-t} \xi_t$. Therefore, by Lemma 32, we have $e^{-ct} M(t) \to \theta W/c$ a.s. By combining this result with Lemma 33, we get

$$\limsup_{t \to \infty} e^{-ct} K(t) \leq \frac{\theta W}{c} - \sum_{k=1}^{N} \frac{\theta W}{c} P(\xi_\tau = k)$$

$$= \frac{\theta W}{c} P(\xi_\tau > N).$$

Letting $N \to \infty$ gives (39). □

REMARK 35. Another consequence of this result is that the proportions of families of size $k$ both in the infinite sites and the infinite alleles models

satisfy

$$\frac{M_k^{\text{GW}}(t)}{M(t)} \to P(\xi_\tau = k),$$

$$\frac{N_k^{\text{GW}}(t)}{M(t)} \to P(\xi_\tau = k),$$

almost surely. We will use this below.

6.4. *Almost sure result for the Beta-coalescent tree.* Let $u > 0$ and consider the reduced tree $\mathbf{T}^u$ at level $u$, which, we recall, has, at level $0 \le t \le 1$, as many vertices as there are excursions between $tu$ and $u$. Suppose mutation marks fall at intensity $\theta\, dt$ per unit length on this tree and for $k \ge 1$, we let $M_k^{\mathbf{T}^u}(t)$ be the number of families of size $k$ at level $0 \le t \le 1$ in the infinite sites model, and let $N_k^{\mathbf{T}^u}(t)$ be the equivalent quantity in the infinite alleles model.

LEMMA 36. *For fixed $u$, conditionally on $\sup_{0 \le s \le T_r} H_s > u$, almost surely as $t \to 0$,*

$$(40) \qquad t^c M_k^{\mathbf{T}^u}(1-t) \to \frac{\theta K}{c} u^{-1/(\alpha-1)} Z_u P(\xi_\tau = k)$$

*and*

$$(41) \qquad t^c N_k^{\mathbf{T}^u}(1-t) \to \frac{\theta K}{c} u^{-1/(\alpha-1)} Z_u P(\xi_\tau = k),$$

*where $K = (\alpha - 1)^{-1/(\alpha-1)}$.*

PROOF. The proof follows from Lemma 34 in exactly the same way that Lemmas 23 and 25 follow from the Kesten–Stigum theorem, the idea being simply that we can again identify $W$ with $K u^{-1/(\alpha-1)} Z_u$ when we look at the reduced tree at level $u$, $\mathbf{T}^u$. □

PROOF OF THEOREM 9. We first note that (40) may be strengthened into the same result where the convergence holds almost surely *for all $\theta$ simultaneously.* Indeed, if we assume that mutation marks come with a label $\theta$ in $(0, \infty)$ and that mutation marks fall on the tree with intensity $d\theta \otimes dt$, where $dt$ stands for the unit length of the tree, we obtain a construction of $M_k^{\mathbf{T}^u}(t)$ for all $\theta$ simultaneously by considering those marks whose label is smaller than $\theta$. (We note for later purposes that, independent of the shape of the tree, such mutation marks may themselves be obtained from a probability measure $Q$ which is a countable collection of independent Poisson processes with intensity $d\theta \otimes dt$.) Observe that since $M_k^{\mathbf{T}^u}(t)$ is monotone in $\theta$ and since (40) holds for all rational $\theta$, it also holds for nonrational values of $\theta$. To get (40) simultaneously for all $\theta$ in the infinite allele case as well, note that $|M_k^{\text{GW}}(t) - N_k^{\text{GW}}(t)| \leq K(t)$ for all $k$ and $t$. Since $K(t)$ is monotone in $\theta$, the result (39) holds for all $\theta > 0$, so (41) also holds simultaneously for all $\theta$.

Let $A_u$ be the event that (40) and (41) hold almost surely for all $\theta$ simultaneously. By applying Lemma 26 with the product probability $P \times Q$, we may assume, without loss of generality, that (40) and (41) hold almost surely for all $\theta$ also at level $u = R^{-1}(1)$, that is,

$$P \times Q(A_{R^{-1}(1)}) = 1.$$

Let $\mathcal{T}_0 = \mathbf{T}^{R^{-1}(1)}$ be the reduced tree at level $R^{-1}(1)$. In order to translate the result to the Beta-coalescent tree, one more fact is needed, since the coalescent tree is not exactly $\mathcal{T}_0$, but a time-change of $\mathcal{T}_0$. (Indeed, for $t \leq 1$, the coalescent tree $\mathcal{T}$ has $t \,\#\mathrm{exc}_{R^{-1}(1-t),R^{-1}(1)}$ branches at level, rather than $\#\mathrm{exc}_{tR^{-1}(1),R^{-1}(1)}$.) In fact, this simply translates into a change of the intensity of the mutation marks for $\mathcal{T}_0$. Indeed, for a given segment in the coalescent tree, between level $R^{-1}(1-t)$ and $R^{-1}(1-s)$ for $s \leq t$, there is a Poisson number of marks with intensity $\theta(t-s)$. So, if $0 \leq \sigma \leq \tau \leq 1$, the number of marks on a segment of the reduced tree $\mathcal{T}_0$ between levels $\sigma$ and $\tau$ is also a Poisson random variable with parameter $\theta(t-s)$, with $R^{-1}(1-t) = \tau R^{-1}(1)$ and $R^{-1}(1-s) = \sigma R^{-1}(1)$. Now, observe that as $t \to 0$ or $\tau \to 1$, this means that the intensity of the marks becomes asymptotic to $\theta R^{-1}(1)/q$, where $q$ is the derivative of the function $R^{-1}(1-t)$ at $t=0$, which was shown to be

$$q = \frac{1}{\alpha(\alpha-1)\Gamma(\alpha)} Z_{R^{-1}(1)}^{\alpha-1}$$

in Lemma 27. Let $M_k^\Pi(t)$ be the number of families of size $k$ obtained from the coalescent tree considered for all $s \geq t$. (I.e., this tree at level $s \geq 0$ has $|\Pi_{t+s}|$ branches.) Using monotonicity of $M_k^{\mathcal{T}_0}(t)$ (number of families of size $k$ in the infinite-site case on $\mathcal{T}_0$) with respect to the intensity, this means that for all $\varepsilon > 0$, for $t$ sufficiently small, $M_k^\Pi(t) \leq M_k^{\mathcal{T}_0}(1 - tq/R^{-1}(1))$, where the intensity is $(\theta + \varepsilon)R^{-1}(1)/q$. Using (40) and the notation $u = R^{-1}(1)$, we have

$$\begin{aligned}\limsup_{t\to 0} t^c M_k^\Pi(t) &\leq \limsup_{t\to 0} t^c M_k^{\mathcal{T}_0}(1 - tq/R^{-1}(1)) \\ &\leq \frac{K(\theta+\varepsilon)u}{qc} u^{-1/(\alpha-1)} Z_u u^c q^{-c} P(\xi_\tau = k) \\ &\leq \frac{\theta+\varepsilon}{c}(\alpha\Gamma(\alpha))^{1/(\alpha-1)} P(\xi_\tau = k)\end{aligned}$$

after simplification, recalling that $c = (2-\alpha)/(\alpha-1)$. We may proceed similarly with the liminf, so we have proven that almost surely as $t \to 0$,

$$t^c M_k^\Pi(t) \to \frac{\theta}{c}(\alpha\Gamma(\alpha))^{1/(\alpha-1)} P(\xi_\tau = k).$$

Combining this result with (39), we get

$$t^c N_k^\Pi(t) \to \frac{\theta}{c}(\alpha\Gamma(\alpha))^{1/(\alpha-1)} P(\xi_\tau = k).$$

The same calculations apply to show that the total number of marks $M^\Pi(t)$ satisfies

$$t^c M^\Pi(t) \to \frac{\theta}{c}(\alpha\Gamma(\alpha))^{1/(\alpha-1)}$$

almost surely. We apply this convergence at times $t = T_n = \inf\{t > 0 : |\Pi_t| \le n\}$. Recall that $T_n \sim (\alpha\Gamma(\alpha))n^{1-\alpha}$ almost surely and that when $|\Pi(T_n)| = n$ (i.e., if the coalescent ever has $n$ blocks), then $M^\Pi(T_n)$ is identical to $M(n)$. On the other hand, by Theorem 1.8 in [8], $\lim_{n\to\infty} P(|\Pi_{T_n}| = n) = \alpha - 1 > 0$, so conditioning on this event which has asymptotically positive probability, we find that

$$n^{\alpha-2} M(n) \to_p \theta \frac{\alpha(\alpha-1)\Gamma(\alpha)}{2-\alpha}$$

(this argument is similar to the one for Theorem 1.9 in [8]). On the other hand, the total number of families $M(n)$ is Poisson with parameter $\theta L_n$, conditionally on $L_n$, where $L_n$ is the total length of the tree, so this gives another proof of Theorem 1.9 in [8] which states that

$$n^{\alpha-2} L_n \to_p \frac{\alpha(\alpha-1)\Gamma(\alpha)}{2-\alpha}.$$

We conclude similarly that

$$n^{\alpha-2} M_k(n) \to_p \theta \frac{\alpha(\alpha-1)\Gamma(\alpha)}{2-\alpha} P(\xi_\tau = k).$$

It follows immediately that the same convergence holds for $N_k(n)$ and this concludes the proof of Theorem 9. □

COROLLARY 37. *Let $K(n)$ be the size of a family chosen uniformly at random among all $M(n)$ families when the population has n individuals. Then*

$$K(n) \to_d \xi_\tau.$$

This is just a reformulation of the fact that the proportions of families of size $k$ converge to $P(\xi_\tau = k)$. Note that $\xi_\tau < \infty$ almost surely, meaning that, asymptotically, a typical family stays of finite size.

**Acknowledgments.** The authors thank Thomas Duquesne and Jean-François Le Gall for helpful discussions, and Jim Pitman for pointing out the connection with [26] and showing us a draft version of this work. N. Berestycki and J. Schweinsberg would also like to express their gratitude to the organizers of the Oberwolfach meeting on Mathematical Population Genetics in August 2005 which was the starting point of the study of the site and allele frequency spectra.

J. BERESTYCKI
LABORATOIRE D'ANALYSE, TOPOLOGIE, PROBABILITÉS
UMR 6632
CENTRE DE MATHÉMATIQUES ET INFORMATIQUE
UNIVERSITÉ DE PROVENCE
39 RUE F. JOLIOT-CURIE
13453 MARSEILLE CEDEX 13
FRANCE
E-MAIL: jberest@cmi.univ-mrs.fr

N. BERESTYCKI
ROOM 121-1984
UNIVERSITY OF BRITISH COLUMBIA
MATHEMATICS ROAD
VANCOUVER, BRITISH COLUMBIA
CANADA V6T 1Z2
E-MAIL: nberestycki@math.ubc.ca

J. SCHWEINSBERG
DEPARTMENT OF MATHEMATICS
U.C. SAN DIEGO
9500 GILMAN DRIVE
LA JOLLA, CALIFORNIA 92093-0112
USA
E-MAIL: jschwein@math.ucsd.edu